\documentclass[12pt]{article}
\usepackage{amssymb,amsmath,theorem}

 \numberwithin{equation}{section}

 \allowdisplaybreaks

\theorempreskipamount=\medskipamount \theorempostskipamount=\medskipamount

\newtheorem{theorem}{Theorem}[section]
\newtheorem{lemma}[theorem]{Lemma}
\newtheorem{corollary}[theorem]{Corollary}

 \let\cal=\mathcal                       
 \let\le\leqslant  \let\ge\geqslant                         
 \let\eps=\varepsilon \let\kappa=\varkappa                  
 \newcommand{\E}{\mathsf{E}\kern 0.07em}                    
 \newcommand{\D}{\mathsf D\kern 0.07em}                     
 \let\emptyset=\varnothing                                  
 \newcommand\qed{\ifhmode\unskip\nobreak\fi\quad            
    \ifmmode\square\else\hbox{$\square$}\fi}                
 \newcommand\proofskip{\vspace{
    \theorempostskipamount}}                                
 \newcommand\pin{\kern.0833em}                              

 \newcommand{\mes}{\mathrm{mes}\kern 0.04em}
 \DeclareMathOperator\dist{dist}
 \newcommand\barV{\kern0.15em\overline{\kern-0.15em V\kern-0.15em}\kern0.15em}

\begin{document}
\vspace*{-3mm}

\section*{The McMillan theorem for colored branching\\[-2pt] processes and dimensions of random fractals}

\medskip

\centerline{\Large\it Victor I. Bakhtin}

\bigskip

\centerline{\large bakhtin@tut.by}

\medskip


\renewcommand{\abstractname}{}
\begin{abstract}
For simplest colored branching processes we prove an analog to the McMillan theorem and calculate
Hausdorff dimensions of random fractals defined in terms of the limit behavior of empirical
measures generated by finite genetic lines. In this setting the role of Shannon's entropy is played
by the Kullback--Leibler divergence and the Hausdorff dimensions are computed by means of the
so-called Billingsley--Kullback entropy, defined in the paper.

\bigskip

{\bf Keywords:} {\it random fractal, Hausdorff dimension, colored brunching process, basin of the
empirical measure, spectral potential, Billingsley--Kullback entropy, Kullback action, maximal
dimension principle}

\medskip

{\bf 2010 Mathematics Subject Classification:\,} 28A80, 37F35, 60J80

\end{abstract}



\bigskip

Let us consider the finite set $X =\{1,\dots,r\}$, whose elements denote different colors, and a
vector $(\mu(1),\dots,\mu(r)) \in [0,1]$. A simplest colored branching process can be defined as an
evolution of a population in which all individuals live the same fixed time and then, when the
lifetime ends, each individual generates (independently of others) a random set of ``children''
containing individuals of colors $1$, \dots, $r$ with probabilities $\mu(1)$, \dots, $\mu(r)$
respectively. We will suppose that the evolution starts with a unique initial individual. It is
suitable to represent this process as a random genealogical tree with individuals as vertices and
each vertex connected by edges with its children. Denote by $X_n$ the set of all genetic lines of
length $n$ (that survive up to generation $n$). The colored branching process can degenerate (when
it turns out that starting from some $n$ all the sets $X_n$ are empty) or, otherwise, evolve
endlessly. Every genetic line $x =(x_1,\dots,x_n)\in X_n$ generates an empirical measure
$\delta_{x,n}$ on the set of colors $X$ by the following rule: for each $i\in X$ the value of
$\delta_{x,n}(i)$ is the fraction of those coordinates of the vector $(x_1,\dots,x_n)$ that
coincide with $i$.

Let $\nu$ be an arbitrary probability measure on $X$. The analog to the McMillan theorem that will
be proved below asserts that under condition of nondegeneracy of the colored branching process the
cardinality of the set $\{\pin x\in X_n\mid \delta_{x,n}\approx\nu\pin\}$ has an almost sure
asymptotics of order $e^{-n\rho(\nu,\mu)}$, where
\begin{equation*}
 \rho(\nu,\mu) =\sum_{i\in X} \nu(i)\ln\frac{\nu(i)}{\mu(i)}.
\end{equation*}
Formally, the value of $\rho(\nu,\mu)$ coincides with the usual Kullback--Leibler divergence and
differs from the latter only in the fact that in our setting the measure $\mu$ is not probability
and so $\rho(\nu,\mu)$ can be negative.

In the paper we investigate also random fractals defined in terms of the sequence of empirical
measures $\delta_{x,n}$ limit behavior. Let $X_\infty$ be the set of infinite genetic lines. Fix an
arbitrary vector $\theta =(\theta(1),\dots,\theta(r)) \in (0,1)^r$ and define the following metrics
on $X_\infty$:
\begin{equation*}
 \dist(x,y) =\prod_{t=1}^n \theta(x_t),\quad\ \text{where}\ \
 n=\inf\pin\{\pin t\mid x_t\ne y_t\pin\} -1.
\end{equation*}

\smallskip

Denote by $V$ any set of probability measures on $X$. It will be proved, in particular, that under
the condition of nondegeneracy of the colored branching process
\begin{equation*}
 \dim_H\pin\{\pin x\in X_\infty \mid \delta_{x,n}\to V\pin\} \pin=\pin
 \sup_{\nu\in V} d(\nu,\mu,\theta)
\end{equation*}
almost surely, where $d(\nu,\mu,\theta)$ is the \emph{Billingsley--Kullback entropy} defined below.

The paper can be divided into two parts. The first one (sections 1--5) contains known results; some
of them have been modified in a certain way for the convenience of use in what follows. Anyway,
most of them are proved below for the completeness and convenience of the reader. The second part
(sections 6--9) contains new results.

In addition, we note that all the results of the paper can be easily extended to Moran's
self-similar geometric constructions in $\mathbb R^n$, but we will not do that.

\section{The spectral potential}\label{1..}

Let $X$ be an arbitrary finite set. Denote by $B(X)$ the space of all real-valued functions on $X$,
by $M(X)$ the set of all positive measures on $X$, and by $M_1(X)$ the collection of all
probability distributions on $X$.

Every measure $\mu\in M(X)$ determines a linear functional on $B(X)$ of the form
\begin{equation*}
 \mu[f] = \int_X f\,d\mu =\sum_{x\in X} f(x)\mu(x).
\end{equation*}
It is easily seen that this functional is \emph{positive} (i.\,e., takes nonnegative values on
nonnegative functions). If, in addition, the measure $\mu$ is probability then this functional is
normalized (takes the value $1$ on the unit function).

Consider the nonlinear functional
\begin{equation}\label{1,,1}
 \lambda(\varphi,\mu) =\ln \mu[e^\varphi],
\end{equation}
where $\varphi\in B(X)$ and $\mu\in M(X)$.  We will call it the \emph{spectral potential}.
Evidently, it is monotone (if $\varphi\ge \psi$ then $\lambda(\varphi,\mu)\ge \lambda(\psi,\mu)$,
additively homogeneous (that is, $\lambda(\varphi+t,\mu) =\lambda(\varphi,\mu) +t$ for each
constant $t$), and analytic in $\varphi$.

Define a family of probability measures $\mu_\varphi$ on $X$, depending on the functional parameter
$\varphi\in B(X)$, by means of the formula
\begin{equation*}
 \mu_{\varphi}[f] = \frac{\mu[e^\varphi f]}{\mu[e^\varphi]}, \qquad f\in B(X).
\end{equation*}
Evidently, each measure $\mu_\varphi$ is equivalent to $\mu$ and has the density $e^{\varphi -
\lambda(\varphi,\mu)}$ with respect to $\mu$.

Let us compute the first two derivatives of the spectral potential with respect to the argument
$\varphi$. Introduce the notation
\begin{equation*}
 \lambda'(\varphi,\mu)[f] =\frac{d\lambda(\varphi+tf,\mu)}{d\pin t}\biggr|_{t=0}.
\end{equation*}
This is nothing more than the derivative of the spectral potential in the direction $f$ at the
point $\varphi$. An elementary computation shows that
\begin{equation}\label{1,,2}
 \lambda'(\varphi,\mu)[f] =\frac{d\ln \mu\bigl[e^{\varphi+tf}\bigr]}{d\pin t}\biggr|_{t=0} =
 \frac{\mu[e^\varphi f]}{\mu[e^\varphi]} = \mu_{\varphi}[f].
\end{equation}
In other words, the derivative $\lambda'(\varphi,\mu)$ coincides with the probability measure
$\mu_\varphi$. Then put
\begin{equation*}
 \lambda''(\varphi,\mu)[f,g] =
 \frac{\partial^2\lambda(\varphi+tf+sg,\mu)}{\partial s\,\partial\pin t}\biggr|_{s,t=0}
\end{equation*}

 \medskip\noindent
and compute this derivative using just obtained formula \eqref{1,,2}:
\begin{equation*}
 \lambda''(\varphi,\mu)[f,g] =\frac{\partial}{\partial s}
 \biggl(\frac{\mu[e^{\varphi+sg}f]}{\mu[e^{\varphi+sg}]}\pin\biggr)\biggr|_{s=0} =
 \frac{\mu[e^\varphi fg]}{\mu[e^\varphi]} -\frac{\mu[e^\varphi f]\pin
 \mu[e^\varphi g]}{(\mu[e^\varphi])^2} = \mu_{\varphi}[fg] -\mu_{\varphi}[f]\pin \mu_{\varphi}[g].
\end{equation*}

In probability theory the expression $\mu_\varphi[t]$ is usually called the expectation of the
random variable $f$ with respect to the probability distribution $\mu_\varphi$, and the expression
$\mu_{\varphi}[fg] -\mu_{\varphi}[f]\pin \mu_{\varphi}[g]$ is called the covariance of random
variables $f$ and $g$. In particular, the second derivative
\begin{equation*}
 \frac{d^2\lambda(\varphi+tf,\mu)}{d\pin t^2}\biggr|_{t=0} =
 \mu_{\varphi}\bigl[f^2\bigr] -\mu_{\varphi}[f]^2
 =\mu_\varphi\bigl[(f -\mu_{\varphi}[f])^2\bigr]
\end{equation*}
is equal to the variance of the random variable $f$ with respect to the distribution $\mu_\varphi$.
Since the variance is nonnegative it follows that the spectral potential is convex in $\varphi$.

\section{The Kullback action}\label{2..}

Denote by $B^*(X)$ the space of all linear functionals on $B(X)$. Then, obviously,
\begin{equation*}
 M_1(X)\subset M(X)\subset B^*(X).
\end{equation*}

The following functional of two arguments $\nu\in B^*(X)$ and $\mu\in M(X)$ will be called the
\emph{Kullback action}:
\begin{equation} \label{2,,1}
 \rho(\nu,\mu) =
 \begin{cases}
   \mu[\varphi\ln\varphi] =\nu[\ln\varphi], &\text{if \,$\nu\in M_1(X)$ \,and \,$\nu=\varphi \mu$},
   \\[2pt]
   +\infty &\text{in all other cases}.
 \end{cases}
\end{equation}
To be more precise, the ``all other cases'' fit into at least one of the three categories: \,a)
singular w.\,r.\,t. $\mu$ probability measures $\nu$, \,b) nonnormalized functionals $\nu$, and
\,c)~nonpositive functionals $\nu$.

In the literature, as far as I know, this functional have been defined only for probability
measures $\nu$ and $\mu$. Different authors call it differently: the relative entropy, the
deviation function, the Kullback--Leibler information function, the Kullback--Leibler divergence.

When $\nu$ is a probability measure the Kullback action can be defined by the explicit formula
\begin{equation} \label{2,,2}
 \rho(\nu,\mu) =\sum_{x\in X} \nu(x)\ln\frac{\nu(x)}{\mu(x)}.
\end{equation}

 \medskip\noindent
In particular, if $\mu(x) \equiv 1$ then the Kullback action differs only in sign from Shannon's
entropy
\begin{equation} \label{2,,3}
 H(\nu) = -\sum_{x\in X} \nu(x)\ln\nu(x).
\end{equation}

\medskip

In the case of probability measure $\mu$ the Kullback action is nonnegative and vanishes only if
$\nu =\mu$. Indeed, if the functional $\nu$ is not an absolutely continuous with respect to $\mu$
probability measure then $\rho(\nu,\mu) =+\infty$. Otherwise, if $\nu$ is a probability measure of
the form $\nu =\varphi\mu$ then from Jensen's inequality and strong convexity of the function $f(x)
=x\ln x$ it follows that
\begin{equation*}
 \rho(\nu,\mu) =\mu[f(\varphi)] \ge f\bigl(\mu[\varphi]\bigr) =0
\end{equation*}
(so long as $\mu[\varphi] =\nu[1] =1$), and the equality $\rho(\nu,\mu) =0$ holds if and only if
$\varphi$ is constant almost everywhere and, respectively, $\nu$ coincides with $\mu$.

Every measure $\mu\in M(X)$ can be put down in the form $\mu =c\mu_1$, where $c =\mu[1]$ and
$\mu_1\in M_1(X)$. If $\nu\in M_1(X)$ then \eqref{2,,2} implies
\begin{equation} \label{2,,4}
 \rho(\nu,\mu) =\rho(\nu,\mu_1) -\ln c \ge -\ln\mu[1].
\end{equation}
In case $\nu\notin M_1(X)$ this inequality holds all the more since the Kullback action is
infinite.

\begin{theorem}\label{2..1}
The spectral potential and the Kullback action satisfy the Young inequality
\begin{equation}\label{2,,5}
 \rho(\nu,\mu)\ge \nu[\psi] -\lambda(\psi,\mu),
\end{equation}

 \medskip\noindent
that turns into equality if and only if\/ $\nu =\mu_\psi$.
\end{theorem}

\emph{Proof.} If $\rho(\nu,\mu) =+\infty$ then the Young inequality is trivial. If $\rho(\nu,\mu)
<+\infty$ then by the definition of Kullback action the functional $\nu$ is an absolutely
continuous probability measure of the form $\nu =\varphi\mu$, where $\varphi$ is a nonnegative
density. In this case
\begin{equation*}
 \lambda(\psi,\mu) \pin=\pin\ln \mu[e^\psi] \pin\ge\pin \ln\! \intop_{\varphi>0} \! e^\psi\, d\mu
 \pin=\pin  \ln\! \intop_{\varphi>0}\! e^{\psi-\ln\varphi}\,d\nu \pin=\pin
 \ln \nu\bigl[e^{\psi -\ln\varphi}\bigr] \pin\ge\pin \nu[\psi-\ln\varphi] \vspace{-2pt}
\end{equation*}
(at the last step we have used Jensen's inequality and concavity of the logarithm function). Since
$\rho(\nu,\mu) =\nu[\ln \varphi]$, this formula implies inequality \eqref{2,,5}.

Recall that $\mu_\psi =e^{\psi-\lambda(\psi,\mu)}\mu$. So if $\nu =\mu_\psi$ then by definition
\begin{equation*}
 \rho(\nu,\mu) =\nu[\psi-\lambda(\psi,\mu)] =\nu[\psi]-\lambda(\psi,\mu).
\end{equation*}
Vice versa, assume that $\rho(\nu,\mu) =\nu[\psi] -\lambda(\psi,\mu)$. Then subtract from the above
equality the Young inequality $\rho(\nu,\mu)\ge \nu[\varphi] -\lambda(\varphi,\mu)$. We obtain
\begin{equation*}
 \lambda(\varphi,\mu) -\lambda(\psi,\mu) \ge \nu[\varphi-\psi].
\end{equation*}
From this follows that $\nu =\lambda'(\psi,\mu)$. Finally, $\lambda'(\psi,\mu)$ coincides with
$\mu_\psi$. \qed

\begin{theorem}\label{2..2}
The Kullback action\/ $\rho(\nu,\mu)$ is the Legendre transform w. r. t.\/ $\nu$ of the spectral
potential\/$:$
\begin{equation}\label{2,,6}
 \rho(\nu,\mu) \pin= \sup_{\psi\in B(X)}\bigl\{\nu[\psi] -\lambda(\psi,\mu)\bigr\}, \qquad
 \nu\in B^*(X),\ \ \mu\in M(X).
\end{equation}
\end{theorem}

\emph{Proof.} By the Young inequality the left hand side of \eqref{2,,6} is not less than the right
one. Therefore it is enough to associate with any functional $\nu\in B^*(X)$ a family of functions
$\psi_t$, depending on the real-valued parameter $t$, on which the equality in \eqref{2,,6} is
attained.

At first, suppose that $\nu$ is an absolutely continuous with respect to $\mu$ probability measure
of the form $\nu =\varphi\mu$, where $\varphi$ is a nonnegative density. Consider the family of
functions
\begin{equation*}
 \psi_t(x) \,=\,
   \begin{cases}
        \ln\varphi(x),&\text{если}\ \ \varphi(x)> 0,\\[2pt]
        -t,&\text{если}\ \ \varphi(x) =0.
   \end{cases}
\end{equation*}

 \medskip\noindent
When $t\to +\infty$ we have the following relations
\begin{gather*}
 \mu\bigl[e^{\psi_t}\bigr] \pin=\intop_{\varphi>0} \! \varphi\,d\mu\pin
    +\intop_{\varphi=0} \!\pin e^{-t}\,d\mu\pin\ \longrightarrow\
    \int_X \varphi\,d\mu \,=\,1,\\[6pt]
 \nu[\psi_t] \pin=\intop_{\varphi>0} \!\varphi\ln\varphi\,d\mu\pin
 +\intop_{\varphi=0} \! -t\varphi\,d\mu \pin=\pin \mu[\varphi\ln\varphi],\\[6pt]
 \nu[\psi_t] -\lambda(\psi_t,\mu) \pin=\pin \nu[\psi_t] -\ln \mu\bigl[e^{\psi_t}\bigr]\,
 \longrightarrow\, \mu[\varphi\ln\varphi] \pin=\pin \rho(\nu,\mu),
\end{gather*}
and so \eqref{2,,6} is proved.

In all the other cases, when $\nu$ is not an absolutely continuous probability measure, by
definition $\rho(\nu,\mu) =+\infty$. Let us examine this cases one after another.

If $\nu$ is a singular relative to $\mu$ probability measure, then there exists $x_0\in X$ such
that $\mu(x_0) =0$ and $\nu(x_0) >0$. In this case consider the family of functions
\begin{equation*}
 \psi_t(x) \,=\,
   \begin{cases}
      t, &\text{если}\ \  x=x_0, \\[2pt]
      0, &\text{если}\ \  x\ne x_0.
   \end{cases}
\end{equation*}
It is easily seen that
\begin{equation*}
 \nu[\psi_t] -\lambda(\psi_t,\mu) \pin\ge\pin t\nu(x_0) -\ln\mu\bigl[e^{\psi_t}\bigr]
 \pin\ge\pin t\nu(x_0) - \ln\mu[1].
\end{equation*}
The right hand side of the above formula goes to $+\infty$ while $t$ increases and \eqref{2,,6}
holds again.

If the functional $\nu$ is not normalized then put $\psi_t =t$. Then the expression
\begin{equation*}
 \nu[\psi_t] -\lambda(\psi_t,\mu)\pin=\pin\nu[t] -\ln\mu[e^t]\pin=\pin t\pin(\nu[1]-1)-\ln\mu[1]
\end{equation*}
is unbounded from the above and hence \eqref{2,,6} is still valid.

Finally, if the functional $\nu$ is not positive then there exists a nonnegative function $\varphi$
such that $\nu[\varphi] <0$. Consider the family $\psi_t =-t\varphi$, where $t>0$. For it
\begin{equation*}
 \nu[\psi_t] -\lambda(\psi_t,\mu) \pin\ge\pin -t\nu[\varphi] -\lambda(0,\mu)
  \, \longrightarrow\, +\infty
\end{equation*}
as $t\to +\infty$, and \eqref{2,,6} remains in force. \qed

\begin{corollary}\label{2..3}
The functional\/ $\rho(\,\cdot\,,\mu)$ is convex and lower semicontinuous on\/ $B^*(X)$.
\end{corollary}

\emph{Proof.} These are properties of the Legendre transform. \qed

\section{The local large deviations principle and\\[-2pt] the McMillan theorem}\label{3..}

As above, we keep to the following notation: $X$ is a finite set, $B(X)$ stands for the space of
real-valued functions on $X$, $B^*(X)$ is the space of linear functionals on $B(X)$, $M_1(X)$ is
the set of all probability measures on $X$, and $M(X)$ is the set of all positive measures on $X$.

To each finite sequence $x =(x_1,\dots,x_n)\in X^n$ let us correspond an \emph{empirical measure}
$\delta_{x,n}\in M_1(X)$ which is supported on the set $\{x_1,\dots,x_n\}$ and assigns to every
point $x_i$ the measure $1/n$. The integral of any function $f$ with respect to $\delta_{x,n}$
looks like
\begin{equation*}
 \delta_{x,n}[f] =\frac{f(x_1)+\,\dotsc\,+f(x_n)}{n}.
\end{equation*}

Denote by $\mu^n$ Cartesian power of a measure $\mu\in M(X)$, which is defined on $X^n$.

\begin{theorem}[\hbox spread -4pt {the local large deviations principle}] \label{3..1}
For any measure\/ $\mu\in M(X)$, any functional\/ $\nu\in B^*(X)$, and\/ $\eps>0$ there exists a
neighborhood\/ $O(\nu)$ such that
\begin{equation}\label{3,,1}
 \mu^n\bigl\{x\in X^n\bigm| \delta_{x,n}\in O(\nu)\bigr\} \pin\le\pin e^{-n(\rho(\nu,\mu) -\eps)}.
\end{equation}
On the other hand, for any\/ $\eps>0$ and any neighborhood\/ $O(\nu)$ the following asymptotic
estimate holds\/$:$
\begin{equation}\label{3,,2}
 \mu^n\bigl\{x\in X^n\bigm| \delta_{x,n}\in O(\nu)\bigr\} \pin\ge\pin e^{-n(\rho(\nu,\mu) +\eps)},
 \qquad n\to\infty.
\end{equation}
\end{theorem}

If $\rho(\nu,\mu) =+\infty$, then by the difference $\rho(\nu,\mu)-\eps$ in \eqref{3,,1} we mean an
arbitrary positive number.

In the case of probability measure $\mu$ Theorem \ref{3..1} is a partial case of Varadhan's large
deviations principle (which explicit formulation can be found, e.\,g., in \cite{Deuschel-Stroock}
and \cite{Varadhan}). Therefore, this theorem can be deduced from Varadhan's large deviations
principle by means of mere renormalization of $\mu$. Nevertheless, we will prove it independently
for the purpose of completeness.

\smallskip

\emph{Proof.} By Theorem \ref{2..2} for any $\eps>0$  there exists $\psi\in B(X)$ such that
\begin{equation}\label{3,,3}
 \rho(\nu,\mu) -\eps/2 \pin<\pin \nu[\psi] -\lambda(\psi,\mu).
\end{equation}
Consider the probability measure $\mu_\psi =e^{\psi-\lambda(\psi,\mu)}\mu$. Obviously,
\begin{equation}\label{3,,4}
 \frac{d\mu^n(x)}{d\mu_\psi^n(x)} \pin=\pin \prod_{i=1}^n \frac{d\mu(x_i)}{d\mu_\psi(x_i)}
 \pin=\pin  \prod_{i=1}^n e^{\lambda(\psi,\mu) -\psi(x_i)} \pin=\pin
 e^{n(\lambda(\psi,\mu) -\delta_{x,n}[\psi])}.
\end{equation}
Define a neighborhood of the functional $\nu$ as follows:
\begin{equation*}
 O(\nu) \pin=\pin \bigl\{\pin \delta\in B^*(X)\bigm| \delta[\psi] >\nu[\psi] -\eps/2\pin\bigr\}.
\end{equation*}
Then it follows from \eqref{3,,4} and \eqref{3,,3} that under the condition $\delta_{x,n}\in
O(\nu)$
\begin{equation*}
 \frac{d\mu^n(x)}{d\mu_\psi^n(x)} \pin<\pin
 e^{n(\lambda(\psi,\mu)-\nu[\psi] +\eps/2)} \pin<\pin e^{n(-\rho(\nu,\mu)+\eps)}.
\end{equation*}
Consequently,
\begin{equation*}
 \mu^n\bigl\{x\in X^n\bigm| \delta_{x,n}\in O(\nu)\bigr\} \, =
 \intop_{\delta_{x,n}\in O(\nu)}\hspace{-1 em} d\mu^n(x) \,\le
 \intop_{\delta_{x,n}\in O(\nu)}\hspace{-1 em} e^{n(-\rho(\nu,\mu)+\eps)}\,d\mu_\psi^n(x)
 \,\le\, e^{-n(\rho(\nu,\mu)-\eps)}.
\end{equation*}
Thus the first part of Theorem \ref{3..1} is proved.

The estimate \eqref{3,,2} is trivial if $\rho(\nu,\mu) =+\infty$. So it is enough to prove it only
in the case when $\nu$ is a probability measure of the form $\nu =\varphi\mu$ and the Kullback
action $\rho(\nu,\mu) =\nu[\ln\varphi]$ is finite. Fix any number $\eps>0$ and neighborhood
$O(\nu)$. Define the sets
\begin{equation*}
 Y_n \pin=\pin \bigl\{\pin x\in X^n\bigm| \delta_{x,n}\in O(\nu),\ \ \big|\delta_{x,n}[\ln\varphi]
 -\nu[\ln\varphi]\big| <\eps/2\pin\bigr\}
\end{equation*}

 \medskip\noindent
(the last inequality in the braces means that $\varphi(x_i)>0$ at each point of the sequence $x
=(x_1,\dots,x_n)$). Note that for $x\in Y_n$
\begin{equation*}
 \frac{d\mu^n(x)}{d\nu^n(x)} \pin=\pin \prod_{i=1}^n \frac{d\mu(x_i)}{d\nu(x_i)} \pin=\pin
 \prod_{i=1}^n \frac{1}{\varphi(x_i)} \pin=\pin e^{-n\delta_{x,n}[\ln\varphi]} \pin>\pin
 e^{-n(\nu[\ln\varphi]+\eps/2)}.
\end{equation*}
Consequently,
\begin{equation}\label{3,,5}
 \mu^n(Y_n) \pin=\pin \int_{Y_n} d\mu^n(x) \pin\ge\pin
 \int_{Y_n} e^{-n(\nu[\ln\varphi]+\eps/2)}\,d\nu^n(x) \pin=\pin
 e^{-n\rho(\nu,\mu) -n\eps/2}\pin\nu^n(Y_n).
\end{equation}
By the Law of large numbers $\nu^n(Y_n)\to 1$. Hence \eqref{3,,5} implies \eqref{3,,2}. \qed

\begin{corollary}[the McMillan theorem] \label{3..2}
For any probability measure\/ $\nu\in M_1(X)$ and\/ $\eps>0$ there exists a neighborhood\/ $O(\nu)$
such that
\begin{equation*}
 \#\{\pin x =(x_1,\dots,x_n)\in X^n\mid \delta_{x,n}\in O(\nu)\pin\} \pin\le\pin e^{n(H(\nu) +\eps)}.
\end{equation*}
On the other hand, for any neighborhood\/ $O(\nu)$ and\/ $\eps>0$
\begin{equation*}
 \#\{\pin x \in X^n\mid \delta_{x,n}\in O(\nu)\pin\} \pin\ge\pin e^{n(H(\nu) -\eps)}
 \quad\ \text{as}\ \ n\to\infty.
\end{equation*}
Here\/ $H(\nu)$ denotes Shannon's entropy defined in\/ \eqref{2,,3}.
\end{corollary}

\emph{Proof.} This follows from equalities \eqref{2,,2}, \eqref{2,,3}, and the previous theorem, if
we set $\mu(x) =1$ for all $x\in X$. \qed

\section{Hausdorff dimension and the maximal dimension principle} \label{4..}

Let us define the Hausdorff dimension of an arbitrary metric space $\Omega$.

Suppose that $\Omega$ is covered by at most countable collection of subsets \pin$\cal U =\{U_i\}$.
Denote by $|\cal U|$ the diameter of this covering: $|\cal U| =\sup |U_i|$, where $|U_i|$ is the
diameter of $U_i$. For every $\alpha\in \mathbb R$ put
\begin{equation*}
 \mes(\cal U,\alpha) =\sum_i |U_i|^\alpha.
\end{equation*}

\medskip

The \emph{Hausdorff measure} (of dimension $\alpha$) of the metric space $\Omega$ is
\begin{equation*}
 \mes(\Omega,\alpha) \pin=\pin \varliminf_{|\kern.04em\cal U|\to 0} \mes(\cal U,\alpha),
\end{equation*}
where \pin$\cal U$ is at most countable covering of $\Omega$. Obviously,
\begin{equation*}
 \mes(\cal U,\beta) \le \mes(\cal U,\alpha)\pin |\cal U|^{\beta-\alpha}
 \quad \text{if}\ \ \beta\ge\alpha.
\end{equation*}
This implies the following property of the Hausdorff measure: if  $\mes(\Omega,\alpha) < \infty$
for some $\alpha$, then $\mes(\Omega,\beta) =0$ for all $\beta> \alpha$.

The \emph{Hausdorff dimension} of the space $\Omega$ is the number
\begin{equation}\label{4,,1}
 \dim_H \Omega \pin=\pin\inf\pin\{\pin \alpha \mid \mes(\Omega,\alpha) =0\pin\}.
\end{equation}
In other words, $\dim_H \Omega =\alpha_0$ if $\mes(\Omega,\alpha) =0$ for all $\alpha>\alpha_0$ and
$\mes(\Omega,\alpha) =\infty$ for all $\alpha<\alpha_0$.

Below we will consider the space of sequences
\begin{equation*}
 X^{\mathbb N} =\{\pin x=(x_1,x_2,x_3,\dots)\pin\}, \quad\ \text{где}\ \ x_i\in X =\{1,\dots,r\}.
\end{equation*}

Let $x =(x_1,x_2,\dots)\in X^{\mathbb N}$. Denote by $Z_n(x)$ the set of sequences
$y=(y_1,y_2,\dots)$ whose first $n$ coordinates coincide with the same coordinates of $x$. This set
will be called a \emph{cylinder of rank} $n$. The collection of all cylinders generates the
\emph{Tychonoff topology} on the space $X^{\mathbb N}$ and the \emph{cylinder $\sigma$-algebra} of
subsets in $X^{\mathbb N}$.

Take an arbitrary positive function $\eta$ on the set of all cylinders that possesses the following
two properties: first, if $Z_n(x) \subset Z_m(y)$ then $\eta(Z_n(x))\le \eta(Z_m(y))$ and, second,
$\eta(Z_n(x)) \to 0$ as $n\to \infty$ at each point $x\in X^{\mathbb N}$. Define the \emph{cylinder
metrics} on $X^{\mathbb N}$ by means of the formula
\begin{equation} \label{4,,2}
 \dist(x,y) =\eta(Z_n(x)), \quad\ \text{where} \ \
 n=\max\pin\{\pin m\mid Z_m(x) =Z_m(y)\pin\}.
\end{equation}
Evidently, the diameter of $Z_n(x)$ in this metrics coincides with $\eta(Z_n(x))$.

Suppose on $X^{\mathbb N}$, besides the cylinder metrics \eqref{4,,2}, a Borel measure $\mu$ is
given. The function
\begin{equation*}
 d_\mu(x) =\varliminf_{n\to\infty} \frac{\ln \mu(Z_n(x))}{\ln |Z_n(x)|}
\end{equation*}

 \medskip\noindent
is called \emph{$($lower\/$)$ pointwise dimension of the measure\/ $\mu$}.

The next theorem provides an effective tool for computing the Hausdorff dimensions of various
subsets of $X^{\mathbb N}$.

\begin{theorem} \label{4..1}
Suppose\/ $A\subset X^{\mathbb N}$. If there exists a finite Borel measure\/ $\mu$ on\/ $X^{\mathbb
N}$ such that\/ $d_\mu(x) \le d$ for each point\/ $x\in A$, then\/ $\dim_H A\le d$. On the
contrary, if\/ $d_\mu(x) \ge d$ for each\/ $x\in A$ and the outer measure\/ $\mu^*(A)$ is positive,
then\/ $\dim_H A \ge d$.
\end{theorem}

It follows that if $d_\mu(x)\equiv d$ on the whole subset $A\subset X^{\mathbb N}$ then its
dimension is equal to $d$.

A weakened version of the second part of Theorem \ref{4..1} in which the condition $d_\mu(x) \ge d$
is replaced by the more strong one $\mu(Z_n(x))\le |Z_n(x)|^d$ is usually called the \emph{mass
distribution principle.}

\smallskip

\emph{Proof.} Every cylinder $Z_n(x)$ is, in fact, a ball in the metrics \eqref{4,,2}, whose radius
equals to its diameter, and vice versa, any ball in this metrics coincides with a cylinder.
Besides, any two cylinders $Z_n(x)$ and $Z_m(y)$ either have empty intersection or one of them is
embedded into other. Therefore, while computing the Hausdorff measure and dimension of a subset
$A\subset X^{\mathbb N}$ it is enough to operate with only disjoint coverings of $A$ by cylinders.

Suppose first that $d_\mu(x) <\alpha$ for all points $x\in A$. Then for each $x\in A$ there exist
arbitrarily small cylinders $Z_n(x)$ satisfying the condition $|Z_n(x)|^\alpha < \mu(Z_n(x))$.
Using this kind of cylinders we can put together a disjoint covering \pin$\cal U$ of the set $A$ of
arbitrarily small diameter. For this covering we have the inequalities
\begin{equation*}
 \mes(\cal U,\alpha) \pin=\sum_{Z_n(x)\in \cal U} |Z_n(x)|^\alpha \pin\le
 \sum_{Z_n(x)\in \cal U} \mu(Z_n(x)) \pin\le\pin \mu\bigl(X^{\mathbb N}\bigr),
\end{equation*}
and hence $\dim_H A\le \alpha$. Thus the first part of the theorem is proved.

Suppose now that $d_\mu(x) >\alpha$ for all points $x\in A$. Define the sets
\begin{equation*}
 A_\eps \pin=\pin \bigl\{ x\in A\bigm| |Z_n(x)|^\alpha > \mu(Z_n(x))\ \ \text{whenever}\ \
 |Z_n(x)| <\eps\bigr\}.
\end{equation*}
Obviously, $A =\bigcup_{\eps>0} A_\eps$. Hence there exists an $\eps$ such that $\mu^*(A_\eps)>0$.
Let \pin$\cal U$ be be a disjoint covering of $A$ by cylinders of diameters less than $\eps$. From
the definition of $A_\eps$ it follows that $\mes(\cal U,\alpha) \ge \mu^*(A_\eps)$. Therefore
$\dim_H A \ge \alpha$, and thus the second part of the theorem is proved. \qed

\medskip

Theorem \ref{4..1} was first proved by Billingsley in the case when the function $\eta$ in
\eqref{4,,2} is a probability measure on $X^{\mathbb N}$ (see \cite[Theorems 2.1 and
2.2]{Billingsley II}). An analog to this theorem for subsets $A\subset \mathbb R^r$ was proved in
\cite{Young} and \cite{Pesin}.

Each point $x =(x_1,x_2,\dots) \in X^{\mathbb N}$ generates a sequence of empirical measures
$\delta_{x,n}$ on the set $X$:
\begin{equation*}
 \delta_{x,n}(i) =\frac{\#\{\pin t\mid x_t=i,\pin\ t\le n\pin\}}{n}, \qquad i\in X.
\end{equation*}

 \medskip\noindent
In other words, $\delta_{x,n}(i)$ is the fraction of those coordinates of the vector
$(x_1,\dots,x_n)$ that coincide with $i$.

For every probability measure $\nu\in M_1(X)$ let us define its \emph{basin} $B(\nu)$ as the set of
all points $x\in X^{\mathbb N}$ such that $\delta_{x,n}$ converges to $\nu$.

Evidently, basins of different measures do not intersect each other and are nonempty. If $x\in
B(\nu)$, and $y\in X^{\mathbb N}$ differs from $x$ in only finite number of coordinates, then $y\in
B(\nu)$. This implies density of each basin in $X^{\mathbb N}$.

Every measure $\nu\in M_1(X)$ generates Bernoulli distribution $P_\nu = \nu^{\mathbb N}$ on the
space $X^{\mathbb N}$. By the strong law of large numbers the basin $B(\nu)$ has probability one
with respect to Bernoulli distribution $P_\nu$, and its complement has zero probability $P_\nu$. In
particular, any basin different from $B(\nu)$ has zero probability.

Points that does not belong to the union of all basins will be called \emph{irregular}. The set or
irregular points has zero probability with respect to any distribution $P_\nu$, where $\nu\in
M_1(X)$. As a result, $X^{\mathbb N}$ turns out to be decomposed into the disjoint union of
different basins and the set of irregular points.

Let us fix some numbers $\theta(i)\in (0,1)$ for all elements $i\in X =\{1,\dots,r\}$, and define a
\emph{cylinder\/ $\theta$-metrics} on $X^{\mathbb N}$ by the rule
\begin{equation}\label{4,,3}
 \dist(x,y) =\prod_{t=1}^n \theta(x_t),\quad\ \text{где}\ \
 n=\inf\pin\{\pin t\mid x_t\ne y_t\pin\} -1.
\end{equation}
It is a partial case of the cylinder metrics \eqref{4,,2}.

For each measure $\nu\in M_1(X)$ and $\theta$-metrics \eqref{4,,3} define the quantity
\begin{equation}\label{4,,4}
 S(\nu,\theta) \pin=\,\frac{\sum_{i=1}^r \nu(i)\ln \nu(i)}{\sum_{i=1}^r \nu(i)\ln \theta(i)}.
\end{equation}
We will call it the \emph{Billingsley entropy} because he was the first who wrote down this formula
and applied it for the computation of Hausdorff dimensions \cite{Billingsley}. He expressed also
this quantity in terms of Shannon's entropy and the Kullback action:
\begin{equation*}
  S(\nu,\theta) \,=\, \frac{H(\nu)}{H(\nu)+\rho(\nu,\theta)}.
\end{equation*}

\begin{theorem} \label{4..2}
Hausdorff dimension of any basin\/ $B(\nu)$ relative to the\/ $\theta$-metrics\/ \eqref{4,,3} is
equal to the Billingsley entropy\/ $S(\nu,\theta)$.
\end{theorem}

A partial case of this theorem in which $\theta(1) =\ldots =\theta(r) =1/r$ was first proved by
Eggleston \cite{Eggleston}. In the complete form this theorem and its generalizations were proved
by Billingsley in \cite{Billingsley,Billingsley II}.

\smallskip

\emph{Proof.} Assume first that $\nu(i)>0$ for every $i=1,\,\dots,\,r$. Obviously,
\begin{equation*}
 \frac{\ln P_\nu(Z_n(x))}{\ln |Z_n(x)|} \pin=\pin
 \frac{\sum_{t=1}^n \ln \nu(x_t)}{\sum_{t=1}^n \ln \theta(x_t)} \pin=\pin
 \frac{\sum_{i=1}^r n\delta_{x,n}(i) \ln \nu(i)}{\sum_{i=1}^r n\delta_{x,n}(i)\ln \theta(i)}.
\end{equation*}
Hence for each point $x\in B(\nu)$ we have
\begin{equation} \label{4,,5}
 d_{P_\nu}(x) \pin=\pin
 \varliminf_{n\to\infty} \frac{\ln P_\nu(Z_n(x))}{\ln |Z_n(x)|}
 \pin=\,\frac{\sum_{i=1}^r \nu(i)\ln \nu(i)}{\sum_{i=1}^r \nu(i)\ln \theta(i)} \pin=\pin S(\nu,\theta).
\end{equation}
Applying Theorem \ref{4..1} to the set $A =B(\nu)$ and measure $\mu =P_\nu$, we obtain the
statement of Theorem \ref{4..2}.

In the general case the same argument provides only lower bound $d_{P_\nu}(x)\ge S(\nu,\theta)$,
that implies the lower bound $\dim_H B(\nu) \ge S(\nu,\theta)$. The inverse inequality is provided
by the next lemma. \qed

\begin{lemma}\label{4..3}
Suppose the space\/ $X^{\mathbb N}$ is equipped with the metrics\/ \eqref{4,,3}. Then for any
measure\/ $\nu\in M_1(X)$ and\/ $\eps>0$ there exists a neighborhood\/ $O(\nu)$ such that Hausdorff
dimension of the set
\begin{equation*}
 A=\bigl\{ x\in X^{\mathbb N}\bigm| \forall\,N\ \exists\,n>N\!:\, \delta_{x,n}\in O(\nu)\bigr\}
\end{equation*}
does not exceed\/ $S(\nu,\theta)+\eps$.
\end{lemma}

\emph{Proof.} Fix a measure $\nu\in M_1(X)$ and an arbitrary positive number $\kappa$. By
McMillan's theorem there exists a neighborhood $O(\nu)$ such that for each positive integer $n$
\begin{equation}\label{4,,6}
 \#\bigl\{ Z_n(x)\bigm| \delta_{x,n}\in O(\nu)\bigr\} \pin\le\pin e^{n(H(\nu) +\kappa)}.
\end{equation}
Decrease this neighborhood in such a way that, in addition, for every measure $\delta\in O(\nu)$
the next inequality holds:
\begin{equation*}
 \sum_{i=1}^r \delta(i)\ln \theta(i) \pin<\pin \sum_{i=1}^r \nu(i) \ln \theta(i) +\kappa.
\end{equation*}
Then for every cylinder $Z_n(x)$ satisfying the condition $\delta_{x,n}\in O(\nu)$ we have the
estimate
\begin{align} \notag
 |Z_n(x)| \pin=\pin \prod_{t=1}^n \theta(x_t) \pin&=\pin\exp\biggl\{\sum_{t=1}^n \ln \theta(x_t)\biggr\}
 \pin=\pin \exp\biggl\{n\sum_{i=1}^r \delta_{x,n}(i)\ln \theta(i)\biggr\} \pin<\pin\\[3pt]
 \pin&<\pin \exp\biggl\{n\sum_{i=1}^r \nu(i)\ln \theta(i) +n\kappa\biggr\}. \label{4,,7}
\end{align}

For any positive integer $N$ the set $A$ is covered by the collection of cylinders
\begin{equation*}
 \cal U_N \pin=\pin\bigcup_{n=N}^\infty \bigl\{ Z_n(x)\bigm| \delta_{x,n}\in O(\nu)\bigr\}.
\end{equation*}
Evidently, the diameter of this covering goes to zero when $N$ increases. Now we can evaluate
$\mes(\cal U_N,\alpha)$ by means of formulas \eqref{4,,6} and \eqref{4,,7}:
\begin{align} \notag
 \mes(\cal U_N,\alpha) \pin&=\sum_{Z_n(x)\in\pin \cal U_N}\hspace{-0.5em} |Z_n(x)|^\alpha
 \pin\le\pin  \sum_{n=N}^\infty e^{n(H(\nu) +\kappa)} \exp\biggl\{\alpha n\sum_{i=1}^r
 \nu(i)\ln \theta(i) +\alpha n\kappa\biggr\} \pin=\pin \\[3pt] \label{4,,8}
 \pin&=\pin \sum_{n=N}^\infty \exp\biggl\{n\pin\biggl(-\sum_{i=1}^r \nu(i)\ln \nu(i)
 +\alpha\sum_{i=1}^r \nu(i)\ln \theta(i) +\kappa +\alpha\kappa\biggr)\!\pin\biggr\}.
\end{align}
If $\alpha > S(\nu,\theta)$, then we can choose so small $\kappa>0$ that the last exponent in
braces is negative, and all the sum \eqref{4,,8} goes to zero as $N\to \infty$. Therefore Hausdorff
measure (of dimension $\alpha$) of the set $A$ is zero, and hence $\dim_H A$ does not exceed
$\alpha$. \qed

\proofskip

We will say that a sequence of empirical measures $\delta_{x,n}$ \emph{condenses} on a subset
$V\subset M_1(X)$ (notation $\delta_{x,n}\succ V$) if it has at least one limit point in $V$.

Similarly to the famous large deviations principle by Varadhan  \cite{Deuschel-Stroock,Varadhan},
it is natural that the next theorem be named the \emph{maximal dimension principle.}

\begin{theorem} \label{4..4}
Let the space\/ $X^{\mathbb N}$ be equipped with the cylinder\/ $\theta$-metrics\/ \eqref{4,,3}.
Then for any nonempty subset\/ $V\subset M_1(X)$
\begin{equation} \label{4,,9}
 \dim_H\pin\bigl\{x\in X^{\mathbb N} \bigm| \delta_{x,n}\succ V\bigr\} \pin=\pin
 \sup_{\nu\in V} S(\nu,\theta).
\end{equation}
\end{theorem}

\emph{Proof.} The set $A =\{\pin x\in X^{\mathbb N} \mid \delta_{x,n}\succ V\pin\}$ contains basins
of all measures $\nu\in V$. So by Theorem \ref{4..2} its dimension is not less than the right hand
side of \eqref{4,,9}.

It is easily seen from the definition \eqref{4,,4} of the Billingsly entropy $S(\nu,\theta)$ that
it depends continuously on the measure $\nu\in M_1(X)$. Consider the closure $\barV$ of $V$.
Obviously, it is compact. Fix any $\eps>0$. By Lemma \ref{4..3} for any measure $\nu\in\barV$ there
exists a neighborhood $O(\nu)$ such that
\begin{equation} \label{4,,10}
 \dim_H\pin\bigl\{x\in X^{\mathbb N}\bigm| \delta_{x,n}\succ O(\nu)\bigr\} \pin\le\pin S(\nu,\theta)+\eps
 \pin\le\pin \sup_{\nu\in V} S(\nu,\theta)+\eps.
\end{equation}
Pick out a finite covering of $\barV$ composed of neighborhoods of this sort. Then the set $A
=\{\pin x\in X^{\mathbb N} \mid \delta_{x,n}\succ V\pin\}$ will be covered by a finite collection
of sets of the form $\{\pin x\in X^{\mathbb N}\mid \delta_{x,n}\succ O(\nu)\pin\}$ satisfying
\eqref{4,,10}. By the arbitrariness of $\eps$ this implies the statement of Theorem \ref{4..4}.
\qed

\proofskip

A very similar to Theorem \ref{4..4} result was proved by Billingsley in \cite[Theorem
7.1]{Billingsley}.

Suppose that a certain subset \pin$\Xi\subset X$ is specified in the set $X =\{1,\dots,r\}$. In
this case the subset \pin$\Xi^{\mathbb N}\subset X^{\mathbb N}$ will be named the \emph{generalized
Cantor set.} It consists of those sequences $x =(x_1,x_2,\dots)$ in which all $x_t\in\Xi$.

\begin{theorem} \label{4..5}
If the space\/ $X^{\mathbb N}$ is equipped with the\/ $\theta$-metrics\/ \eqref{4,,3} then
Hausdorff dimension of the generalized Cantor set\/ \pin$\Xi^{\mathbb N}$ coincides with the unique
solution of Moran's equation
\begin{equation} \label{4,,11}
 \sum_{i\pin\in\pin\Xi} \theta(i)^s =1.
\end{equation}
\end{theorem}

\medskip

This theorem was first proved by Moran in 1946 \cite{Moran} for generalized Cantor subsets of the
real axis and afterwards it was extended by Hutchinson \cite{Hutchinson} to the attractors of
self-similar geometric constructions in $\mathbb R^r$. Let us show how it can be derived from the
maximal dimension principle.

\smallskip

\emph{Proof.} Let $s$ be the solution to Moran's equation. Introduce a probability distribution
$\nu$ on $X$, setting $\nu(i) =\theta(i)^s$ for $i\in\Xi$ and $\nu(i) =0$ for $i\notin\Xi$. Then
\begin{equation} \label{4,,12}
 S(\nu,\theta) \pin=\pin \frac{\sum_{i=1}^r \nu(i)\ln \nu(i)}{\sum_{i=1}^r \nu(i)\ln \theta(i)} \pin=\pin
 \frac{\sum_{i\in\Xi}\theta(i)^s\ln \theta(i)^s}{\sum_{i\in\Xi}\theta(i)^s\ln \theta(i)} \pin=\pin s.
\end{equation}

Consider the set $B(\nu)\cap\Xi^{\mathbb N}$. It has the unit measure with respect to the
distribution $P_\nu =\nu^{\mathbb N}$. Besides, for every point $x\in B(\nu)\cap\Xi^{\mathbb N}$ by
\eqref{4,,5} we have the equality $d_{P_\nu}(x) =S(\nu,\theta)$. In this setting it follows from
Theorem \ref{4..1} and formula \eqref{4,,12} that
\begin{equation} \label{4,,13}
 \dim_H \bigl(B(\nu)\cap\Xi^{\mathbb N}\bigr) =S(\nu,\theta) =s.
\end{equation}

Denote by $V$ the collection of all probability measures on $X$ supported on \pin$\Xi\subset X$.
Evidently, for each point $x\in\Xi^{\mathbb N}$ all the limit points of the sequence $\delta_{x,n}$
belong to $V$. Hence, we can apply Theorem \ref{4..4} that implies
\begin{equation} \label{4,,14}
 \dim_H \pin\Xi^{\mathbb N} \pin\le\pin
 \dim_H\pin\bigl\{x\in X^{\mathbb N} \bigm| \delta_{x,n}\succ V\bigr\} \pin=\pin
 \sup_{\nu\in V} S(\nu,\theta) \pin=\pin
 \sup_{\nu\in V}\frac{\sum_{i\in\Xi} \nu(i)\ln \nu(i)}{\sum_{i\in\Xi} \nu(i)\ln \theta(i)}\pin.
\end{equation}
Note that for every measure $\nu\in V$
\begin{equation*}
 s\sum_{i\pin\in\pin\Xi} \nu(i)\ln \theta(i) \pin-\sum_{i\pin\in\pin\Xi} \nu(i)\ln \nu(i) \pin=\pin
 \sum_{\nu(i)>0}\hspace{-0.3em} \nu(i)\ln\frac{\theta(i)^s}{\nu(i)} \pin\le\pin
 \ln\Biggl\{\sum_{\nu(i)>0}\hspace{-0.3em} \nu(i)\frac{\theta(i)^s}{\nu(i)}\!\pin\Biggr\} \pin\le\pin
 0,
\end{equation*}
where we have used concavity of the logarithm function. It follows that the right hand side in
\eqref{4,,14} does not exceed $s$. Finally, comparing \eqref{4,,13} and \eqref{4,,14}, we obtain
the desired equality $\dim_H \pin\Xi^{\mathbb N} =s$. \qed

\section{Brunching processes} \label{5..}

First let us introduce the basic notions about the simplest Galton--Watson brunching process.

Suppose that a random variable $Z$ takes nonnegative values $k\in\mathbb Z_+$ with probabilities
$p_k$. The \emph{Galton--Watson brunching process} is a sequence of integer-valued random variables
$Z_0$, $Z_1$, $Z_2$, \dots{} such that $Z_0\equiv 1$, \,$Z_1 =Z$, and further each $Z_{n+1}$ is
defined as the sum of $Z_n$ independent counterparts of the random variable $Z$. In particular, if
$Z_n =0$ then $Z_{n+1} =0$ as well. Usually $Z_n$ is thought of as the total number of descendants
in $n$-th generation of a unique common ancestor under the condition that each descendant
independently of others gives birth to $Z$ children.

It is known that in some cases the posterity of the initial ancestor may degenerate (when starting
from a certain $n$ all $Z_n$ are zeros) and in other cases it can ``flourish'' (when $Z_n$ grows
exponentially). The type of behavior of the brunching process depends on the mean number of
children of any individual
\begin{equation*}
 m =\E Z = \sum_{k=0}^\infty kp_k
\end{equation*}
and on the generating function of that number
\begin{equation*}
 f(s) =f_1(s) =\sum_{k=0}^\infty p_k s^k.
\end{equation*}
Obviously, the restriction of the function $f(s)$ to the segment $[0,1]$ is nonnegative,
nondecreasing, convex, and satisfies $f(1) =1$ and $f'(1) =m$.

In the theory of brunching processes (see, for instance, \cite{Athreya-Ney,Harris}) the following
statements have been proved.

\begin{theorem} \label{5..1}
The generating functions of the number of descendants in\/ $n$-th generation
\begin{equation*}
 f_n(s) =\sum_{k=0}^\infty \mathsf{P}\{Z_n=k\}\pin s^k
\end{equation*}

 \medskip\noindent
satisfy the recursion relation\/ $f_{n+1}(s) =f(f_n(s))$.
\end{theorem}

\begin{theorem} \label{5..2}
If\/ $m\le 1$ then the brunching process degenerates almost surely\/ $($except the case when each
individual gives birth to exactly one child\/$)$. If\/ $m>1$ then the probability\/ $q$ of
degeneration is less than\/ $1$ and coincides with a unique nonunit root of the equation\/ $f(s)
=s$ on the segment\/ $[0,1]$.
\end{theorem}

\begin{theorem} \label{5..3}
If\/ $m>1$ and\/ $\E Z^2<\infty$ then the sequence\/ $W_n =Z_n/m^n$ converges almost surely to a
random variable\/ $W$ such that\/ $\mathsf{P}\{W>0\} =1-q$. If\/ $m>1$ and\/ $\E Z^2 =\infty$ then
for any number\/ $m'<m$ with probability\/ $1-q$
\begin{equation*}
 \lim_{n\to\infty} Z_n/m^n <\infty, \qquad \lim_{n\to\infty} Z_n\big/(m')^n =\infty
\end{equation*}
$($here\/ $q$ is the probability of degeneration of the brunching process\/$)$.
\end{theorem}

Thereby, in the case $m>1$ there is an alternative for the total number of descendants $Z_n$:
either it vanishes at a certain moment $n_0$ (with probability $q<1$) or it is asymptotically
equivalent to $W m^n$ (with the complementary probability $1-q$), where the random variable $W>0$
does not depend on $n$ (except the case $\E Z^2 =\infty$, when only the logarithmic equivalence
$\ln Z_n \sim \ln m^n$ is guaranteed). All other types of the descendants' number behavior have
zero probability.

We will exploit these theorems in the study of colored brunching processes.

Suppose now that each individual may give birth to children of $r$ different colors (or $r$
different genders, if one likes). We will suppose that the posterity of each individual in the
first generation represents a random set $X$ containing random number $k_1$ of children of the
first color, random number $k_2$ of children of the second color, and so on up to $k_r$ children of
color $r$. All elements of $X$ (including elements of the same color) are treated as different. The
ordered array $k=(k_1,k_2,\dots,k_r) \in \mathbb Z_+^r$ will be called the \emph{color structure}
of the set of children $X$. Denote by $p_k$ the probability of birth of the set $X$ with color
structure $k=(k_1,k_2,\dots,k_r)$. Naturally, all the probabilities $p_k$ are nonnegative and
\begin{equation*}
 \sum_{k\in\mathbb Z_+^r} p_k =1.
\end{equation*}

 \smallskip\noindent
If an individual $x_1$ gave birth to $x_2$, then $x_2$ gave birth to $x_3$, and so on up to an
individual $x_n$, then the sequence $x =(x_1,\dots,x_n)$ will be called the \emph{genetic line} of
length $n$.

Let us construct a new branching process taking into account not only the total number of
descendants but also the color of each individual and all its upward and downward lineal relations.
This process may be thought of as a random \emph{genealogical tree} with a common ancestor in the
root and all its descendants in the vertices, where each parent is linked with all its children. In
the case of degenerating population its genealogical tree is finite, and in the case of
``flourishing'' one the tree is infinite.

Formally it is convenient to define such a process as a sequence of random sets $X_n$ containing
all genetic lines of length $n$. As the first set $X_1$ we take $X$. The subsequent $X_n$ are built
up by induction: if $X_n$ is already known, then for all genetic lines $(x_1,\dots,x_n)\in X_n$
define disjoint independent random sets of children $X(x_1,\dots,x_n)$, each with color structure
distribution as in $X$, and put
\begin{equation*}
 X_{n+1} \pin=\pin \bigl\{ (x_1,\dots,x_n,x_{n+1})\bigm| (x_1,\dots,x_n)\in X_n,\ \, x_{n+1}\in
 X(x_1,\dots,x_n) \bigr\}.
\end{equation*}
The built in such a way stochastic process $X_1$, $X_2$, $X_3$, \dots\ will be referred to as the
\emph{colored branching process} (or \emph{unconditional} colored branching process if one wishes
to emphasize that the posterity of any individual is independent of its color and genealogy).

\section{The McMillan theorem for colored branching\\[-2pt] processes} \label{6..}

Consider a colored branching process $X_1$, $X_2$, \dots\ determined by a finite collection of
colors $\Omega =\{1,\dots,r\}$ and a probability distribution $\{\pin p_k\mid k\in\mathbb
Z^r_+\pin\}$, where $k=(k_1,\dots,k_r)$ is the color structure of each individual's set of children
$X$. We will always think that $X_1$ is generated by a unique initial individual.

For any genetic line $x =(x_1,\dots,x_n)\in X_n$ define the \emph{spectrum} $\delta_{x,n}$ as the
corresponding empirical measure on $\Omega$ by the rule
\begin{equation} \label{6,,1}
 \delta_{x,n}(i) =\frac{\#\{\pin t\mid g(x_t) =i\pin\}}{n}, \qquad i\in\Omega,
\end{equation}
where $g(x_t)$ denotes the color of $x_t$. In other words, $\delta_{x,n}(i)$ is the fraction of
individuals of color $i$ in the genetic line $x$. Our next goal is to obtain asymptotical estimates
for cardinalities of the random sets
\begin{equation*}
 \bigl\{ x =(x_1,\dots,x_n)\in X_n\bigm| \delta_{x,n}\in O(\nu)\bigr\},
\end{equation*}
where $O(\nu)$ is a small neighborhood of the distribution $\nu$ on the set of colors $\Omega$.

Denote by $\mu(i)$ the expectation of members of color $i$ in $X$:
\begin{equation} \label{6,,2}
 \mu(i) =\sum_{k\in\mathbb Z^r_+} k_ip_k, \qquad i=1,\,\dots,\,r.
\end{equation}
Provided all $\mu(i)$ are finite, the vector $\mu =(\mu(1),\dots,\mu(r))$ can be regarded as a
measure on the set of colors $\Omega$. This measure generates the measure $\mu^n$ on $\Omega^n$ as
Cartesian product.

Define a mapping $G:X_n\to\Omega^n$ by means of the formula
\begin{equation*}
 G(x_1,\dots,x_n) =\bigl(g(x_1),\pin\dots,\pin g(x_n)\bigr),
\end{equation*}
where $g(x_t)$ is the color of $x_t$.

\begin{lemma} \label{6..1}
For any\/ $\omega =(\omega_1,\dots,\omega_n) \in \Omega^n$ we have
\begin{equation} \label{6,,3}
 \E\kern 0.05em\#\{\pin x\in X_n\mid G(x) =\omega\pin\} \pin=\pin \prod_{t=1}^n \mu(\omega_t)
 \pin=\pin \mu^n(\omega).
\end{equation}
\end{lemma}

\emph{Proof.} Cast out the last coordinate in $\omega$ and let $\omega' = (\omega_1, \dots,
\omega_{n-1})$. For any genetic line $(x_1,\dots,x_{n-1})\in X_{n-1}$, by virtue of the definition
of unconditional colored branching process we have
\begin{equation*}
 \E\kern 0.05em\#\{\pin x_{n}\in X(x_1,\dots,x_{n-1})\mid g(x_{n}) =\omega_{n}\pin\}
 \pin=\pin \mu(\omega_{n}).
\end{equation*}
Evidently, this expression does not depend on $x' =(x_1,\dots,x_{n-1})$. Therefore,
\begin{equation*}
 \E\kern 0.05em\#\{\pin x\in X_{n}\mid G(x) =\omega\pin\} \pin=\pin
 \E\kern 0.05em\#\{\pin x'\in X_{n-1}\mid G(x') =\omega'\pin\}\pin \mu(\omega_{n}).
\end{equation*}
Repeated application of the latter equality gives \eqref{6,,3}. \qed

\proofskip

Define for the measure $\mu$ from \eqref{6,,2} the Kullback action
\begin{equation*}
 \rho(\nu,\mu) =\sum_{i\in\Omega} \nu(i)\ln\frac{\nu(i)}{\mu(i)}, \qquad \nu\in M_1(\Omega),
\end{equation*}
where $M_1(\Omega)$ is the set of all probability measures on $\Omega$. This formula is a copy of
\eqref{2,,2}.

\begin{theorem} \label{6..2}
Suppose\/ $X_1,\, X_2,\, \dots$ is an unconditional colored brunching process with finite
collection of colors\/ $\Omega$. Then for any\/ $\eps>0$ and probability measure\/ $\nu\in
M_1(\Omega)$ there exists a neighborhood\/ $O(\nu)\subset M_1(\Omega)$ such that for all natural\/
$n$
\begin{equation} \label{6,,4}
 \E\kern 0.05em\#\{\pin x\in X_n\mid \delta_{x,n}\in O(\nu)\pin\} \pin\le\pin
 e^{n(-\rho(\nu,\mu)+\eps)}.
\end{equation}
On the other hand, for any\/ $\eps>0$ and any neighborhood\/ $O(\nu)$
\begin{equation} \label{6,,5}
 \E\kern 0.05em\#\{\pin x\in X_n\mid \delta_{x,n}\in O(\nu)\pin\} \pin\ge\pin
 e^{n(-\rho(\nu,\mu)-\eps)} \quad\ \text{as}\ \ n\to\infty.
\end{equation}
\end{theorem}

If $\rho(\nu,\mu) =+\infty$, the expression  $-\rho(\nu,\mu)+\eps$ in \eqref{6,,4} should be
treated as an arbitrary negative real number.

\smallskip

\emph{Proof.} It follows from \eqref{6,,1} that for every genetic line $x\in X_n$ its spectrum
$\delta_{x,n}$ coincides with the empirical measure $\delta_{\omega,n}$, where $\omega =G(x)$.
Therefore,
\begin{equation} \label{6,,6}
 \#\{\pin x\in X_n\mid \delta_{x,n}\in O(\nu)\pin\} \,=
 \sum_{\omega\in\Omega^n:\pin \delta_{\omega,n}\in O(\nu)}
 \#\{\pin x\in X_n\mid G(x) =\omega\pin\}. \\[-6pt]
\end{equation}
It follows from \eqref{6,,3} and \eqref{6,,6} that
\begin{equation*}
 \E\kern 0.05em\#\{\pin x\in X_n\mid \delta_{x,n}\in O(\nu)\pin\} \pin=\pin
 \mu^n\{\pin\omega\in\Omega^n\mid \delta_{\omega,n}\in O(\nu)\pin\}.
\end{equation*}
The latter equality converts estimates \eqref{6,,4} and \eqref{6,,5} into already proved estimates
\eqref{3,,1}, \eqref{3,,2} from the large deviations principle. \qed

\proofskip

Remarkable that the last reference to the large deviations principle serves a unique ``umbilical
cord'' linking the first three sections of the paper with others.

Now we are ready to state an analog of the McMillan theorem for colored branching processes. Let
$q^*$ be a probability of degeneration of the process (probability of the occasion that starting
from a certain number $n$ all the sets $X_n$ turn out to be empty).

\begin{theorem} \label{6..3}
Suppose\/ $X_1,\, X_2,\, \dots$ is an unconditional colored brunching process with finite
collection of colors\/ $\Omega$. Then for any\/ $\eps>0$ and any probability measure\/ $\nu\in
M_1(\Omega)$ there exists a neighborhood\/ $O(\nu)\subset M_1(\Omega)$ such that for almost sure
\begin{equation} \label{6,,7}
 \#\{\pin x\in X_n\mid \delta_{x,n}\in O(\nu)\pin\} \pin<\pin e^{n(-\rho(\nu,\mu)+\eps)} \quad\
 \text{as}\ \ n\to\infty.
\end{equation}
On the other hand, if\/ $\rho(\nu,\mu) <0$ then for any neighborhood\/ $O(\nu)$ and positive\/
$\eps$ the estimate
\begin{equation} \label{6,,8}
 \#\{\pin x\in X_n\mid \delta_{x,n}\in O(\nu)\pin\} \pin>\pin e^{n(-\rho(\nu,\mu)-\eps)}
 \quad\ \text{as}\ \ n\to\infty
\end{equation}

 \medskip\noindent
holds with probability\/ $1-q^*$ $($or almost surely under the condition that our branching process
does not degenerate\/$)$.
\end{theorem}

\emph{Proof.} Application of Chebyshev's inequality to \eqref{6,,4} gives
\begin{equation*}
 \mathsf{P}\bigl\{\#\{\pin x\in X_n\mid \delta_{x,n}\in O(\nu)\pin\} \ge
 e^{n(-\rho(\nu,\mu)+2\eps)}\bigr\} \pin\le\pin e^{-n\eps}.
\end{equation*}
Sum up these inequalities over all $n\ge N$:
\begin{equation*}
 \mathsf{P}\bigl\{\exists\, n\ge N\!:\  \#\{\pin x\in X_n\mid \delta_{x,n}\in O(\nu)\pin\} \ge
 e^{n(-\rho(\nu,\mu)+2\eps)}\bigr\} \pin\le\pin \frac{e^{-N\eps}}{1-e^{-\eps}}.
\end{equation*}
This implies \eqref{6,,7} with constant $2\eps$ instead of $\eps$, that does not change its sense.

Proceed to the second part of the theorem. Let $\kappa = -\rho(\nu,\mu) -\eps$ and the number
$\eps$ be so small that $\kappa>0$. By the second part of Theorem \ref{6..2} for any neighborhood
$O(\nu)$ there exists $N$ such that
\begin{equation*}
 \E\kern 0.05em\#\{\pin x\in X_N\mid \delta_{x,N}\in O(\nu)\pin\} \pin>\pin e^{N\kappa}.
\end{equation*}
Without loss of generality we may assume that $O(\nu)$ is convex.

Construct a Galton--Watson branching process satisfying the conditions
\begin{gather} \label{6,,9}
 Z_1 \pin=\pin \#\{\pin x\in X_N\mid \delta_{x,N}\in O(\nu)\pin\}, \\[6pt] \label{6,,10}
 Z_n \pin\le\pin \#\{\pin x\in X_{nN}\mid \delta_{x,nN}\in O(\nu)\pin\},
 \qquad n=2,\,3,\,\dots
\end{gather}
Let the random variable $Z_1$ be defined by \eqref{6,,9}. For $n>1$ define $Z_n$ as a total number
of genetic lines
\begin{equation*}
 (x_1,\pin\dots,\pin x_N,\ \ .\ \ .\ \ .\ \ ,\pin x_{(n-1)N+1},\pin\dots,\pin x_{nN})
 \pin\in\pin X_{nN}
\end{equation*}
such that the spectrum of each segment $(x_{kN+1},\pin\dots,\pin x_{(k+1)N})$ belongs to $O(\nu)$.
In other words, we will treat as ``individuals'' of the process $Z_1$, $Z_2$, $Z_3$, \dots\ those
segments $(x_{kN+1},\pin\dots,\pin x_{(k+1)N})$ of genetic lines of the initial process whose
spectrum lies in $O(\nu)$. Then \eqref{6,,10} follows from convexity of $O(\nu)$, and from
unconditionality of the initial colored branching process it can be concluded that the sequence
$Z_1$, $Z_2$, \dots\ in fact forms a Galton--Watson branching process.

By construction, $\E Z_1 >e^{N\kappa}$. In this setting Theorem \ref{5..3} asserts that there is an
alternative for the sequence $Z_n$: either it tends to zero with a certain probability $q<1$ or it
grows faster than $e^{nN\kappa}$ with probability $1-q$. In the second case, by virtue of
\eqref{6,,10},
\begin{equation} \label{6,,11}
 \#\{\pin x\in X_{nN}\mid \delta_{x,nN}\in O(\nu)\pin\}\, e^{-nN\kappa}
 \pin\rightarrow\pin \infty \quad\ \text{при}\ \ n\to\infty.
\end{equation}

\medskip

To finish the proof we have to do two things: verify that in fact \eqref{6,,11} is valid with
probability $1-q^*$ and get rid of the multiplier $N$ there. To do this we will exploit two ideas.
First, if the colored branching process $X_1$, $X_2$, \dots\ were generated by $m$ initial
individuals instead of the unique one, then \eqref{6,,11} would be valid with probability at least
$1-q^m$. Second, if one genetic line is a part of another and the ratio of their lengthes is close
to $1$ then their spectra are close as well.

Obviously, the total number of individuals in the $n$-th generation of the initial branching
process $X_1$, $X_2$, $X_3$, \dots\ equals $|X_n|$. The sequence of random variables $|X_n|$ forms
a Galton--Watson brunching process with probability of degeneration $q^*$, that does not exceed
$q$. Therefore, the sequence $|X_n|$ grows exponentially with probability $1-q^*$.

Consider the colored brunching process $X_{k+1}$, $X_{k+2}$, $X_{k+3}$, \dots\ obtained from the
initial one by virtue of truncation of the first $k$ generations. It represents a union of $|X_k|$
independent brunching processes generated by all individuals of $k$-th generation. It satisfies
\eqref{6,,11} with probability at least $1- q^{|X_k|}$. Hence for the initial process with even
greater probability we obtain the condition
\begin{equation} \label{6,,12}
 \#\{\pin x\in X_{k+nN}\mid \delta_{x,k+nN}\in O^*(\nu)\pin\}\, e^{-nN\kappa}
 \pin\rightarrow\pin \infty \quad\ \text{при}\ \ n\to\infty,
\end{equation}
where $O^*(\nu)$ is an arbitrary neighborhood of $\nu$ containing the closure of $O(\nu)$.

Suppose the sequence $|X_n|$ grows exponentially. Then for every $m\in\mathbb N$ define the numbers
\begin{equation*}
 k_i \pin=\pin\min\pin\{\pin k: |X_k|\ge m,\ \ k=i \hspace{-0.6em}\mod N\pin\},
 \qquad i=0,\,1,\,\dots,\,N-1.
\end{equation*}
For each $k=k_i$, the condition \eqref{6,,12} holds with probability at least $1-q^m$, and in
common they give the estimate
\begin{equation*}
 \#\{\pin x\in X_n\mid \delta_{x,n}\in O^*(\nu)\pin\} \pin>\pin e^{n\kappa}
 \quad\ \text{as}\ \ n\to\infty
\end{equation*}
with probability at least $1-Nq^m$. By virtue of the arbitrariness of $m$ this estimate is valid
almost surely (under the condition $|X_n|\to\infty$, which takes place with probability $1-q^*$).
It is equivalent to \eqref{6,,8}. \qed

\section{Dimensions of random fractals (upper bounds)}\label{7..}

We proceed investigation of the colored brunching process $X_1$, $X_2$, \dots\ with finite
collection of colors $\Omega =\{1,\dots,r\}$. Let us consider the corresponding set of infinite
genetic lines
\begin{equation*}
 X_\infty \pin=\pin\bigl\{ x=(x_1,x_2,x_3,\dots)\bigm| (x_1,\dots,x_n)\in X_n\ \
 \forall\,n\in \mathbb N \bigr\}.
\end{equation*}

Define the cylinder $\theta$-metrics on $X_\infty$
\begin{equation}\label{7,,1}
 \dist(x,y) =\prod_{t=1}^n\theta(x_t), \qquad n=\inf\pin\{\pin t\mid x_t\ne y_t\pin\} -1,
\end{equation}
where the numbers $\theta(1)$, \dots, $\theta(r)$ are taken from $(0,1)$.

We will be interested in Hausdorff dimensions of both the space $X_\infty$ and its various subsets
defined in terms of partial limits of empirical measures on $\Omega$ (those measures are called
spectra and denoted $\delta_{x,n}$). If the colored brunching process degenerates then $X_\infty$
is empty. Therefore of interest is only the case when $m =\E |X_1| >1$ and  the cardinality of
$X_n$ increases with rate of order $m^n$.

As before, denote by $\mu(i)$, where $i\in \Omega$, the expectation of individuals of color $i$ in
the random set $X_1$. It will be always supposed that $\mu(i)<\infty$. Consider any probability
measure $\nu\in M_1(\Omega)$. It will be proved below that the dimension of $\{\pin x\in
X_\infty\mid \delta_{x,n}\to \nu\pin\}$ can be computed by means of the function
\begin{equation} \label{7,,2}
 d(\nu,\mu,\theta) \pin=\pin \frac{\rho(\nu,\mu)\strut}{\sum_{i=1}^r \nu(i)\ln \theta(i)\strut}
 \pin=\pin \frac{\sum_{i=1}^r \nu(i)\ln\frac{\textstyle \nu(i)}{\textstyle \mu(i)}}{\sum_{i=1}^r
 \nu(i)\ln \theta(i)\strut}.
\end{equation}

 \vspace{1pt}\noindent
We will name it the \emph{Billingsley--Kullback entropy.}

In \eqref{7,,2} the numerator is the Kullback action and the denominator is negative. If $\mu$ is a
probability measure on $\Omega$ then the Kullback action is nonnegative. But in our setting this is
not the case since $m =\mu(1)+\ldots+ \mu(r) > 1$. In particular, if $\mu(i)>\nu(i)$ for all $i\in
\Omega$ then the Kullback action will be negative, and the Billingsley--Kullback entropy positive.
Note, in addition, that if $\mu(1) =\ldots =\mu(r) =1$ then $-\rho(\nu,\mu)$ is equal to Shannon's
entropy $H(\nu)$, and the whole of Billingsley--Kullback entropy turns into the Billingsley entropy
\eqref{4,,4}.

\begin{lemma}\label{7..1}
Let the space\/ $X_\infty$ of infinite genetic lines be equipped with the metrics\/ \eqref{7,,1}.
Then for any probability measure\/ $\nu\in M_1(\Omega)$ and any\/ $\eps>0$ there exists a
neighborhood\/ $O(\nu)$ such that Hausdorff dimension of the set
\begin{equation*}
 A =\bigl\{ x\in X_\infty\bigm| \forall\,N\ \exists\,n>N\!:\, \delta_{x,n}\in O(\nu)\bigr\}
\end{equation*}
does not exceed\/ $d(\nu,\mu,\theta)+\eps$ almost surely.
\end{lemma}

\emph{Proof} is carried out in the same manner as in Lemma \ref{4..3}. Take any $\kappa>0$. By
Theorem \ref{6,,4} there exists a neighborhood  $O(\nu)$ such that almost surely
\begin{equation}\label{7,,3}
 \#\bigl\{ x\in X_n\bigm| \delta_{x,n}\in O(\nu)\bigr\} \pin\le\pin e^{n(-\rho(\nu,\mu) +\kappa)}
 \quad\ \text{as}\ \ n\to\infty.
\end{equation}
Reduce this neighborhood in such a way that in addition for all measures $\delta\in O(\nu)$,
\begin{equation*}
 \sum_{i=1}^r \delta(i)\ln \theta(i) \pin<\pin \sum_{i=1}^r \nu(i) \ln \theta(i) +\kappa.
\end{equation*}
Then for each cylinder $Z_n(x)$ satisfying the condition $\delta_{x,n}\in O(\nu)$ we have the
estimate
\begin{align} \notag
 |Z_n(x)| \pin=\pin \prod_{t=1}^n \theta(x_t) \pin&=\pin\exp\biggl\{\sum_{t=1}^n \ln \theta(x_t)\biggr\}
 \pin=\pin \exp\biggl\{n\sum_{i=1}^r \delta_{x,n}(i)\ln \theta(i)\biggr\} \pin<\pin\\[3pt]
 \pin&<\pin \exp\biggl\{n\sum_{i=1}^r \nu(i)\ln \theta(i) +n\kappa\biggr\}. \label{7,,4}
\end{align}

For every natural $N$ the set $A$ is covered by the collection of cylinders
\begin{equation*}
 \cal U_N \pin=\pin\bigcup_{n=N}^\infty \bigl\{ Z_n(x)\bigm| \delta_{x,n}\in O(\nu)\bigr\}.
\end{equation*}
Evidently, the diameter of this covering tends to zero as $N\to \infty$. Hence $\mes(\cal
U_N,\alpha)$ can be estimated by virtue of formulas \eqref{7,,3} and \eqref{7,,4}:
\begin{align} \notag
 \mes(\cal U_N,\alpha) \pin&=\sum_{Z_n(x)\in\pin \cal U_N}\hspace{-0.5em} |Z_n(x)|^\alpha
 \pin\le\pin  \sum_{n=N}^\infty e^{n(-\rho(\nu,\mu) +\kappa)}
 \exp\biggl\{\alpha n\sum_{i=1}^r \nu(i)\ln \theta(i) +\alpha n\kappa\biggr\} \pin=\pin \\[3pt]
 \pin&=\pin \sum_{n=N}^\infty \exp\biggl\{n\pin\biggl(-\sum_{i=1}^r \nu(i)\ln\frac{\nu(i)}{\mu(i)}
 +\alpha\sum_{i=1}^r \nu(i)\ln \theta(i)+\kappa+\alpha\kappa\biggr)\!\pin\biggr\}. \label{7,,5}
\end{align}
If $\alpha > d(\nu,\mu,\theta)$ then $\kappa$ can be chosen so small that the last exponent in
braces is negative, and all the sum \eqref{7,,5} tends to zero as $N\to \infty$. Therefore
Hausdorff measure (of dimension $\alpha$) of the set $A$ is zero, and its dimension does not exceed
$\alpha$. \qed

\proofskip

As before, we say that the sequence of empirical measures $\delta_{x,n}$ condenses on a subset
$V\subset M_1(\Omega)$ (notation $\delta_{x,n}\succ V$) if it has a limit point in $V$.

\begin{theorem} \label{7..2}
Let\/ $X_1,\, X_2,\, X_3,\, \dots$ be an unconditional colored brunching process with finite set of
colors\/ $\Omega$, and the set\/ $X_\infty$  of all infinite genetic lines equipped with the
cylinder metrics\/ \eqref{7,,1}. Then for any subset\/ $V\subset M_1(\Omega)$ almost surely
\begin{equation}\label{7,,6}
  \dim_H\pin\{\pin x\in X_\infty \mid \delta_{x,n}\succ V\pin\} \pin\le\pin
  \sup_{\nu\in V} d(\nu,\mu,\theta).
\end{equation}
In particular, $\dim_H X_\infty\le s$ for almost sure, where\/ $s$ is a unique root of the ``Bowen
equation''
\begin{equation} \label{7,,7}
 \sum_{i=1}^r \mu(i) \theta(i)^s =1.
\end{equation}
\end{theorem}

\medskip

\emph{Proof.} It follows from the definition of the Billingsley--Kullback entropy
$d(\nu,\mu,\theta)$ that it depends continuously on the measure $\nu\in M_1(\Omega)$. Let $\barV$
be the closure of $V$. Obviously, it is compact. Take an arbitrary $\eps>0$. By Lemma \ref{7..1}
for any measure $\nu\in\barV$ there exists a neighborhood $O(\nu)$ such that almost surely
\begin{equation} \label{7,,8}
 \dim_H\pin\bigl\{x\in X_\infty\bigm| \delta_{x,n}\succ O(\nu)\bigr\} \pin\le\pin
 d(\nu,\mu,\theta)+\eps \pin\le\pin \sup_{\nu\in V} d(\nu,\mu,\theta)+\eps.
\end{equation}
Choose a finite covering of $\barV$ by neighborhoods of this kind. Then the set $\{\pin x\in
X_\infty \mid \delta_{x,n}\succ V\pin\}$ will be covered by a finite collection of sets of the form
$\{\pin x\in X_\infty \mid \delta_{x,n}\succ O(\nu)\pin\}$ satisfying \eqref{7,,8}. By the
arbitrariness of $\eps$ this implies the first statement of Theorem \ref{7..2}.

Let $s$ be a solution of equation \eqref{7,,7}. Note that for any measure $\nu\in M_1(\Omega)$,
since the logarithm function is concave,
\begin{gather*}
 s\sum_{i=1}^r \nu(i)\ln \theta(i) \pin-\sum_{i=1}^r \nu(i)\ln\frac{\nu(i)}{\mu(i)} \pin=\pin
 \sum_{i=1}^r \nu(i)\ln\frac{\mu(i)\theta(i)^{s}}{\nu(i)} \pin\le\pin \\[6pt]
 \pin\le\pin  \ln\Biggl\{\sum_{\nu(i)>0}\hspace{-0.3em}
 \nu(i)\frac{\mu(i)\theta(i)^{s}}{\nu(i)}\!\pin\Biggr\} \pin\le\pin 0.
\end{gather*}
Consequently, $d(\nu,\mu,\theta)\le s$. Now the second part of our theorem follows from the first
one if we take $V=M_1(\Omega)$. \qed

\proofskip

{\bf Remark.\,} In fact the ``Bowen equation'' is an equation of the form $P(s\varphi) =0$, where
$P(s\varphi)$ is the topological pressure of a weight function $s\varphi$ in a dynamical system
(more detailed explanations can be found in \cite{Pesin}). If we replace the topological pressure
$P(s\varphi)$ by the spectral potential
\begin{equation*}
 \lambda(s\varphi,\mu) = \ln\sum_{i=1}^r e^{s\varphi(i)}\mu(i),
 \quad\ \text{where}\ \ \varphi(i) =\ln \theta(i),
\end{equation*}
then the Bowen equation turns into the equation $\lambda(s\varphi,\mu) =0$, which is equivalent to
\eqref{7,,7}.

\section{Block selections of colored brunching processes} \label{8..}

Let $\xi_1$, $\xi_2$, $\xi_3$, \dots\ be a sequence if independent identically distributed random
variables taking values $0$ or $1$ (independent Bernoulli trials).

\begin{lemma} \label{8..1}
If\/ $0<p'<p<1$ and\/ $\mathsf{P}\{\xi_i=1\}\ge p$, then
\begin{equation} \label{8,,1}
 \mathsf{P}\{\pin\xi_1+\ldots+\xi_k \ge p'k\pin\} \pin\to\pin 1 \quad \text{as}\ \ k\to\infty
\end{equation}
uniformly with respect to the probability\/ $\mathsf{P}\{\xi_i=1\}\ge p$.
\end{lemma}

\emph{Proof.} In the case $\mathsf{P}\{\xi_i=1\} =p$ this follows from the law of large numbers. If
$\mathsf{P}\{\xi_i=1\}$ increases then the probability in the left hand side of \eqref{8,,1}
increases as well. \qed

\proofskip

Consider a colored brunching process $X_1$, $X_2$, \dots\ with finite set of colors $\Omega =
\{1,\dots,r\}$. Each $X_n$ consists of genetic lines $(x_1,x_2,\dots,x_n)$ of length $n$, in which
every subsequent individual has been born by the previous. Fix a (large enough) natural $N$. We
will split genetic lines of length divisible by $N$ into blocks of length $N$:
\begin{equation*}
 (x_1,x_2,\dots,x_{nN}) =(y_1,\dots,y_n), \quad\ \text{where}\quad
 y_{k} =(x_{(k-1)N+1},\pin\dots,\pin x_{kN}).
\end{equation*}
Each block $y_k$ generates an empirical measure $\delta_{y_k}$ (spectrum) on $\Omega$ by the rule
\begin{equation*}
 \delta_{y_k}(i) =\frac{\#\{\pin t\mid g(x_t)=i,\ \, (k-1)N <t\le kN\pin\}}{N},
\end{equation*}
where $g(x_t)$ denotes the color of $x_t$.

A \emph{block selection of order\/ $N$} from a colored brunching process $X_1$, $X_2$, \dots\ is
any sequence of random subsets $Y_n\subset X_{nN}$ with the following property: if
$(y_1,\dots,y_{n+1})\in Y_{n+1}$ then $(y_1,\dots,y_n)\in Y_n$. In this case the sequence of blocks
$(y_1,\dots,y_{n+1})$ will be called a \emph{prolongation} of the sequence $(y_1,\dots,y_n)$.

As above (see \eqref{6,,2}), denote by $\mu(i)$ the expectation of children of color $i$ born by
each individual, and by $\mu$ the corresponding measure on $\Omega$.

\begin{theorem} \label{8..2}
Let\/ $X_1,\, X_2,\, X_3,\, \dots$ be an unconditional colored brunching process with finite set of
colors\/ $\Omega$ and probability of degeneration\/ $q^*<1$. If a measure\/ $\nu\in M_1(\Omega)$
satisfies the condition\/ $\rho(\nu,\mu) <0$, then for any its neighborhood\/ $O(\nu)\subset
M_1(\Omega)$ and any number\/ $\eps>0$ with probability\/ $1-q^*$ one can extract from the
brunching process a block selection\/ $Y_1$, $Y_2$, \dots\ of an order\/ $N$ such that each
sequence of blocks\/ $(y_1,\dots,y_n)\in Y_n$ has at least\/ $l(N)$ prolongations in\/ $Y_{n+1}$,
where
\begin{equation} \label{8,,2}
 l(N) = e^{N(-\rho(\nu,\mu)-\eps)},
\end{equation}
and the spectra of all blocks belong to\/ $O(\nu)$.
\end{theorem}

\emph{Proof.} Fix any numbers $p$ and $\eps$ satisfying the conditions
\begin{equation*}
 0<p<p+\eps<1-q^*, \qquad \rho(\nu,\mu)+\eps<0.
\end{equation*}
By the second part of Theorem \ref{6..3} for all large enough $N$ we have
\begin{equation} \label{8,,3}
 \mathsf{P}\bigl\{ \#\{\pin x\in X_N \mid \delta_{x,N}\in O(\nu)\pin\} >
 e^{N(-\rho(\nu,\mu)-\eps/2)} \bigr\} \pin\ge\pin p +\eps.
\end{equation}

Further we will consider finite sequences of random sets $X_1$, \dots, $X_{nN}$ and extract from
them block selections $Y_1$, \dots, $Y_n$ of order $N$ such that the spectra of all their blocks
belong to $O(\nu)$ and each sequence of blocks $(y_1,\dots,y_k)\in Y_k$ has at least $l(N)$
prolongations in $Y_{k+1}$. Denote by $A_n$ the event of existence of a block selection with these
properties. Define one more event $A$ by the condition
\begin{equation*}
 \#\{\pin x\in X_N \mid \delta_{x,N}\in O(\nu)\pin\} \pin>\pin l(N)\pin e^{N\eps/2}.
\end{equation*}

It follows from \eqref{8,,2} and \eqref{8,,3} that $\mathsf{P}(A)\ge p+\eps$. Evidently, $A\subset
A_1$. Therefore, $\mathsf{P}(A_1)\ge p+\eps$. Now we are going to prove by induction that
$\mathsf{P}(A_n)\ge p$ whenever the order $N$ of selection is large enough. Let us perform the step
of induction. Assume that $\mathsf{P}(A_n)\ge p$ is valid for some $n$. Consider the conditional
probability $\mathsf{P}(A_{n+1}|A)$. By the definition of events $A_{n+1}$ and $A$ it cannot be
less than the probability of the following event: there are at least $l(N)$ wins in a sequence of
$[l(N) e^{N\eps/2}]$ independent Bernoulli trials with probability of win $\mathsf{P}(A_n)$ in
each. Using Lemma \ref{8..1} (with $p' =p/2$ and $k=[l(N) e^{N\eps/2}]$) one can make this
probability greater than $1-\eps$ at the expense of increasing~$N$. Then,
\begin{equation*}
 \mathsf{P}(A_{n+1}) \ge \mathsf{P}(A)\pin\mathsf{P}(A_{n+1}|A) > (p+\eps)(1-\eps) >p.
\end{equation*}
Thus the inequality $\mathsf{P}(A_n) >p$ is proved for all $n$.

It means that with probability greater than $p$ one can extract from the sequence $X_1$, \dots,
$X_{nN}$ a block selection $Y_1$, \dots, $Y_n$ of order $N$ such that the spectra of all blocks
belong to the neighborhood $O(\nu)$ and each sequence of blocks $(y_1,\dots,y_k)\in Y_k$ has at
least $l(N)$ prolongations in $Y_{k+1}$.

To obtain a block selection of infinite length with the same properties, we will construct finite
block selections $Y_1$, \dots, $Y_n$ in the following manner. Initially, suppose that every $Y_k$,
where $k\le n$, consists of all sequences of blocks $(y_1,\dots,y_k)\in X_{kN}$ such that the
spectrum of each block lies in $O(\nu)$. At the first step we exclude from $Y_{n-1}$ all sequences
of blocks having less than $l(N)$ prolongations in $Y_n$, and then exclude from $Y_n$ all
prolongations of the sequences that have been excluded from $Y_{n-1}$. At the second step we
exclude from $Y_{n-2}$ all sequences of blocks having after the first step less than $l(N)$
prolongations in the modified $Y_{n-1}$, and then exclude from $Y_{n-1}$ and $Y_{n}$ all
prolongations of the sequences that have been excluded from $Y_{n-2}$. Proceeding further in the
same manner, after $n$ steps we will obtain a block selection $Y_1$, \dots, $Y_n$ such that each
sequence of blocks from any $Y_k$ has at least $l(N)$ prolongations in $Y_{k+1}$. Evidently, this
selection will be the maximal among all selections of order $N$ having the mentioned property.
Therefore with probability at least $p$ all the sets $Y_k$ are nonempty.

For every $n$ let us construct, as is described above, the maximal block selection $Y^{(n)}_1$,
\dots, $Y^{(n)}_n$. From the maximality of these selections it follows that
\begin{equation*}
 Y^{(n)}_n\supset Y^{(n+1)}_n\supset Y^{(n+2)}_n\supset\ldots
\end{equation*}
Define the sets $Y_n =\bigcap_{k\ge n} Y^{(k)}_n$. Then with probability at least $p$ all of them
are nonempty and compose an infinite block selection from Theorem \ref{8..2}. Since $p$ may be
chosen arbitrarily close to $1-q^*$, such selections do exist with probability $1-q^*$. \qed

\proofskip

Theorem \ref{8..2} can be strengthened by taking several measures in place of a unique measure
$\nu\in M_1(\Omega)$.

\begin{theorem} \label{8..3}
Let\/ $X_1,\, X_2,\, X_3,\, \dots$ be an unconditional colored brunching process with finite set of
colors\/ $\Omega$ and probability of degeneration\/ $q^*<1$. If a finite collection of measures\/
$\nu_i\in M_1(\Omega)$, where\/ $i=1,\,\dots,\,k$, satisfy the inequalities\/ $\rho(\nu_i,\mu) <0$,
then for any neighborhoods\/ $O(\nu_i)\subset M_1(\Omega)$ and any\/ $\eps>0$ with probability\/
$1-q^*$ one can extract from the brunching process a block selection\/ $Y_1$, $Y_2$, \dots\ of an
order\/ $N$ such that for every\/ $i=1,\,\dots,\,k$ each sequence of blocks\/ $(y_1,\dots,y_n)\in
Y_n$ has at least\/
\begin{equation*}
 l_i(N) = e^{N(-\rho(\nu_i,\mu)-\eps)}
\end{equation*}
prolongations\/ $(y_1,\dots,y_n,y)\in Y_{n+1}$ with the property\/ $\delta_{y}\in O(\nu_i)$.
\end{theorem}

It can be proved in the same manner as the previous one, only now the event $A_n$ should be
understood as existence of a finite block selection $Y_1$, \dots, $Y_n$ satisfying the conclusion
of Theorem \ref{8..3} and the event $A$ should be defined by the system of inequalities
\begin{equation*}
 \#\{\pin x\in X_N \mid \delta_{x,N}\in O(\nu_i)\pin\} \pin>\pin l_i(N)\pin e^{N\eps/2},
 \qquad i=1,\,\dots,\,k.
\end{equation*}
We leave details to the reader.

\section{Dimensions of random fractals (lower bounds)}\label{9..}

Now we proceed investigation of the space of infinite genetic lines
\begin{equation*}
 X_\infty \pin=\pin\bigl\{ x=(x_1,x_2,x_3,\dots)\bigm| (x_1,\dots,x_n)\in X_n\ \
 \forall\,n\in \mathbb N \bigr\},
\end{equation*}
which is generated by an unconditional colored brunching process $X_1$, $X_2$, \dots\ with finite
set of colors $\Omega =\{1,\dots,r\}$. It is supposed that there is a measure
\begin{equation*}
 \mu = (\mu(1),\,\dots,\,\mu(r))
\end{equation*}
on $\Omega$, where $\mu(i)$ denotes the expectation of children of color $i$ born by each
individual, and $X_\infty$ is equipped with the cylinder $\theta$-metrics \eqref{7,,1}.

\begin{theorem} \label{9..1}
Let\/ $X_1,\, X_2,\, X_3,\, \dots$ be an unconditional colored brunching process with finite set of
colors\/ $\Omega$ and probability of degeneration\/ $q^*<1$. If a measure\/ $\nu\in M_1(\Omega)$
satisfies the condition\/ $d(\nu,\mu,\theta) >0$, then with probability\/ $1-q^*$ for any
neighborhood\/ $O(\nu)$ we have the lower bound
\begin{equation} \label{9,,1}
 \dim_H\pin \{\pin x\in X_\infty\mid \exists\,N\ \forall\,n>N\ \, \delta_{x,n}\in O(\nu)\pin\}
 \,\ge\, d(\nu,\mu,\theta).
\end{equation}
\end{theorem}

\emph{Proof.} Fix any number $\alpha <d(\nu,\mu,\theta)$ and so small $\eps>0$ that
\begin{equation} \label{9,,2}
 d(\nu,\mu,\theta) \pin=\pin \frac{\rho(\nu,\mu)\strut}{\sum_{i=1}^r \nu(i)\ln \theta(i)\strut} \pin>\pin
 \frac{\rho(\nu,\mu) +2\eps\strut}{\sum_{i=1}^r \nu(i)\ln \theta(i) -\eps\strut} \pin>\pin \alpha.
\end{equation}
Then choose a convex neighborhood $O^*(\nu)$ whose closure lies in $O(\nu)$ such that for any
measure $\delta\in O^*(\nu)$
\begin{equation} \label{9,,3}
 \sum_{i=1}^r \delta(i)\ln \theta(i) >\sum_{i=1}^r \nu(i)\ln \theta(i)-\eps.
\end{equation}

 \smallskip

By Theorem \ref{8..2} with probability $1-q^*$ one can extract from the brunching process under
consideration a block selection $Y_1$, $Y_2$, \dots\ of order $N$ such that any sequence of blocks
$(y_1,\dots,y_n)\in Y_n$ has at least $l(N)$ prolongations in $Y_{n+1}$, where
\begin{equation*}
 l(N) = e^{N(-\rho(\nu,\mu)-\eps)},
\end{equation*}
and for each block $y_k =(x_{(k-1)N+1},\dots,x_{kN})$ the corresponding empirical measure
$\delta_{y_k}$ (spectrum) belongs to $O^*(\nu)$. Exclude from this selection a certain part of
genetic lines in such a way that each of the remaining sequences of blocks $(y_1,\dots,y_n)\in Y_n$
would have exactly $[l(N)]$ prolongations in $Y_{n+1}$.

Define the random set
\begin{equation*}
 Y_\infty \pin=\pin \bigl\{ y=(y_1,y_2,\dots) \bigm| (y_1,\dots,y_n)\in Y_n, \ \,
 n=1,\,2,\,\dots\bigr\}.
\end{equation*}
Any sequence $y=(y_1,y_2,\dots)\in Y_\infty$ consists of blocks of length $N$. Having written down
in order the elements of all these blocks, we obtain from $y$ an infinite genetic line $x =
(x_1,x_2,\dots)\in X_\infty$. Denote it as $\pi(y)$. By the definition of $Y_\infty$ the spectrum
of each block $y_k$ belongs to $O^*(\nu)$. For every point $x=\pi(y)$, where $y\in Y_\infty$, the
empirical measure $\delta_{x,nN}$ is an arithmetical mean of empirical measures corresponding to
the first $n$ blocks of $y$, and so belongs to $O^*(\nu)$ as well. It follows that
\begin{equation} \label{9,,4}
 \pi(Y_\infty) \pin\subset\pin \{\pin x\in X_\infty\mid \exists\,N\ \forall\,n>N\ \,
 \delta_{x,n}\in O(\nu)\pin\}.
\end{equation}

The family of all cylinders of the form $Z_{nN}(x)$, where $x\in \pi(Y_\infty)$, generates some
$\sigma$-algebra on $\pi(Y_\infty)$. Define a probability measure $P$ on this $\sigma$-algebra such
that
\begin{equation*}
 P\bigl(Z_{nN}(x)\bigr) =[l(N)]^{-n}.
\end{equation*}
Then for all large enough $N$, all $x \in \pi(Y_\infty)$, and all natural $n$
\begin{equation*}
 P\bigl(Z_{nN}(x)\bigr) \le e^{nN(\rho(\nu,\mu)+2\eps)}.
\end{equation*}
On the other hand, by \eqref{9,,3}
\begin{gather*}
 |Z_{nN}(x)| \pin=\pin \prod_{t=1}^{nN} \theta(x_t) \pin=\pin
 \exp\biggl\{\sum_{t=1}^{nN} \ln \theta(x_t)\biggr\} \pin=\pin
 \exp\biggl\{\sum_{i=1}^r nN\delta_{x,nN}(i)\ln \theta(i)\biggr\} \pin\ge\pin \\[3pt]
 \pin\ge\pin \exp\biggl\{nN\biggl(\pin\sum_{i=1}^r \nu(i)\ln \theta(i) -\eps\biggr)\biggr\}.
\end{gather*}
It follows from the last two formulas and \eqref{9,,2} that
\begin{equation*}
 |Z_{nN}(x)|^\alpha \pin\ge\pin
 \exp\biggl\{nN\alpha\biggl(\pin\sum_{i=1}^r \nu(i)\ln \theta(i)-\eps\biggr)\biggr\} \pin\ge\pin
 e^{nN(\rho(\nu,\mu)+2\eps)} \pin\ge\pin P\bigl(Z_{nN}(x)\bigr).
\end{equation*}

Now we are ready to compute the Hausdorff measure of dimension $\alpha$ of the set $\pi(Y_\infty)$.
If, while computing the Hausdorff measure, we used coverings of $\pi(Y_\infty)$ not with any
cylinders, but with only cylinders of orders divisible by $N$, then the last formula would imply
that such a measure will be at least $P(\pi(Y_\infty)) =1$. Any cylinder can be put in a cylinder
of order divisible by $N$ such that the difference of their orders will be less than $N$ and the
ratio of their diameters greater than $\min \theta(i)^{N}$. Therefore,
\begin{equation*}
 \mes\bigl(\pi(Y_\infty),\alpha\bigr) \ge \min \theta(i)^{N\alpha}
\end{equation*}
and hence $\dim_H \pi(Y_\infty) \ge \alpha$.

The set defined in the right hand part of \eqref{9,,4} contains $\pi(Y_\infty)$. Then its dimension
is at least $\alpha$ too. Recall that we have proved this fact by means of a block selection that
exists with probability $1-q^*$. By the arbitrariness of $\alpha<d(\nu,\mu,\theta)$ this implies
the desired bound \eqref{9,,1} with the same probability. \qed

\begin{theorem} \label{9..2}
Let\/ $s$ be a root of the Bowen equation
\begin{equation*}
 \sum_{i=1}^r \mu(i) \theta(i)^s =1.
\end{equation*}
If\/ $s\le 0$, then\/ $X_\infty =\emptyset$ almost surely. Otherwise, if\/ $s>0$, then\/ $X_\infty$
is nonempty with a positive probability, and with the same probability its dimension equals\/ $s$.
\end{theorem}

\emph{Proof.} The expectation of total number of children of each individual in the brunching
process generating the set $X_\infty$ is equal to $m=\mu(1)+\,\ldots\,+\mu(r)$. If $s\le 0$, then
$m\le 1$. In this case by Theorem \ref{5..2} our brunching process degenerates almost surely, and
$X_\infty =\emptyset$.

If $s>0$, then $m>1$. In this case by Theorem ref{5..2} our brunching process is degenerate with a
positive probability, and $X_\infty$ is nonempty. Define a measure $\nu\in M_1(\Omega)$ by means of
the equality
\begin{equation*}
 \nu(i) =\mu(i)\theta(i)^{s}, \quad\  i\in \Omega.
\end{equation*}
Then, evidently, $d(\nu,\mu,\theta) =s$. By the previous Theorem $\dim_H X_\infty \ge s$ with the
same probability with which $X_\infty \ne\emptyset$. On the other hand, by Theorem \ref{7..2} the
inverse inequality holds almost surely. \qed

\proofskip

A more general version of Theorem \ref{9..2}, in which the similarity coefficients $\theta(1)$,
\dots, $\theta(r)$ are random, is proved in
\cite{Falconer-article,Falconer-book,Graf,Mauldin-Williams}.

For every probability measure $\nu\in M_1(\Omega)$ define a basin $B(\nu)\subset X_\infty$ as the
set of all infinite genetic lines $x =(x_1,x_2,x_3,\dots)$ such that the corresponding sequence of
empirical measures $\delta_{x,n}$ converges to $\nu$. What is the dimension of $B(\nu)$? By Theorem
\ref{7..2} it does not exceed the Billingsley--Kullback entropy $d(\nu,\mu,\theta)$ with
probability $1$. On the other hand, the inverse inequality does not follow from the previous
results (and, in particular, from Theorem \ref{9..1}). To obtain it, we ought to enhance the
machinery of block selections.

\begin{lemma} \label{9..3}
Let\/ $Q_1$, \dots, $Q_{2^r}$ be vertices of a cube in\/ $\mathbb R^r$. Then there exists a choice
law\/ $i\!:\mathbb R^r\to \{1,\dots,2^r\}$ such that if neighborhoods\/ $O(Q_i)$ are small enough,
sequences
\begin{equation*}
 \delta_n\in\mathbb R^r \quad \text{and}\quad
 \Delta_n =\frac{\delta_1+\ldots+\delta_n}{n}
\end{equation*}
satisfy the conditions\/ $\delta_{n+1}\in O\bigl(Q_{i(\Delta_n)}\bigr)$ and\/ $\delta_1\in O(Q_1)
\cup \ldots\cup O(Q_{2^r})$, then the sequence\/ $\Delta_n$ converges to the center of the cube.
\end{lemma}

\emph{Proof.} First consider the case $r=1$, when the cube turns to a segment. Let, for
definiteness, $Q_1 =-1$ and $Q_2 =1$. Set
\begin{equation} \label{9,,5}
 i(\Delta) =
   \begin{cases}
      1,& \text{если}\ \ \Delta\ge 0,\\[1pt]
      2,& \text{если}\ \ \Delta<0.
   \end{cases}
\end{equation}
Take any neighborhoods $O(Q_1)$ and $O(Q_2)$ with radii at most $1$. Then for any sequence
$\delta_n$ satisfying the conditions $\delta_{n+1}\in O\bigl(Q_{i(\Delta_n)}\bigr)$ and
$|\delta_1|<2$ we have the estimate $|\Delta_n|<2/n$. It may be easily proved by induction. Thus in
the one-dimensional case the lemma is proved. To prove it in the multidimensional case one should
choose a coordinate system with origin at the center of the cube and axes parallel to edges of the
cube and apply the choice law \eqref{9,,5} to each of the coordinates independently. \qed

\begin{theorem} \label{9..4}
Let\/ $X_1,\, X_2,\, X_3,\, \dots$ be an unconditional colored brunching process with finite set of
colors\/ $\Omega$ and probability of degeneration\/ $q^*<1$. If a measure\/ $\nu\in M_1(\Omega)$
satisfies the condition\/ $d(\nu,\mu,\theta) >0$, then with probability\/ $1-q^*$
\begin{equation*}
 \dim_H B(\nu) = d(\nu,\mu,\theta).
\end{equation*}
\end{theorem}

\emph{Proof.} Fix any number $\alpha <d(\nu,\mu,\theta)$ and so small $\eps>0$ that
\begin{equation*}
 d(\nu,\mu,\theta) \pin=\pin \frac{\rho(\nu,\mu)\strut}{\sum_{i=1}^r \nu(i)\ln \theta(i)\strut} \pin>\pin
 \frac{\rho(\nu,\mu) +3\eps\strut}{\sum_{i=1}^r \nu(i)\ln \theta(i) -\eps\strut} \pin>\pin \alpha.
\end{equation*}

The set $M_1(\Omega)$ is in fact a simplex of dimension $r-1$, where $r =|\Omega|$. Suppose first
that $\nu$ is an inner point of this simplex (in other words, $\nu(i)>0$ for all $i\in \Omega$).
Take a small convex neighborhood $O(\nu)\subset M_1(\Omega)$ such that for any measure $\delta\in
O(\nu)$
\begin{gather*}
 \rho(\delta,\mu) < \rho(\nu,\mu) +\eps,\\[3pt]
 \sum_{i=1}^r \delta(i)\ln \theta(i) >\sum_{i=1}^r \nu(i)\ln \theta(i)-\eps.
\end{gather*}

Let $Q_1$, \dots, $Q_{2^{r-1}}$ be vertices of some cube in $O(\nu)$ with center at $\nu$. Define
for them small neighborhoods $O(Q_i)\subset O(\nu)$ as in Lemma \ref{9..3}. By Theorem \ref{8..3}
with probability $1-q^*$ one can extract from the colored branching process $X_1$, $X_2$, \dots\ a
block selection $Y_1$, $Y_2$, \dots\ of order $N$ such that for every $i\le 2^{r-1}$ each sequence
of blocks $(y_1,\dots,y_n)\in Y_n$ has at least
\begin{equation*}
 l(N) = e^{N(-\rho(\nu,\mu)-2\eps)}
\end{equation*}
prolongations $(y_1,\dots,y_n,y)\in Y_{n+1}$ possessing the property $\delta_{y}\in O(Q_i)$.

Exclude from this block selection a certain part of genetic lines so that each of the remaining
sequences of blocks $(y_1,\dots,y_n)\in Y_n$ would have exactly $[l(N)]$ prolongations
$(y_1,\dots,y_n,y)\in Y_{n+1}$, and all these prolongations would satisfy the choice law from Lemma
\ref{9..3}, namely,
\begin{equation*}
 \delta_y\in O\bigl(Q_{i(\Delta_n)}\bigr), \quad\ \text{where} \quad
 \Delta_n =\frac{\delta_{y_1}+\ldots+\delta_{y_n}}{n}.
\end{equation*}

Denote by $\pi(Y_\infty)$ the set of all infinite genetic lines $(x_1,x_2,\dots)\in X_\infty$ for
which every initial segment of length $nN$, has been partitioned into blocks of length $N$, turns
into an element of $Y_n$. Then by Lemma \ref{9..3} we have the inclusion $\pi(Y_\infty) \subset
B(\nu)$.

Reproducing reasoning from the proof of Theorem \ref{9..1} one can ascertain that the dimension of
$\pi(Y_\infty)$ is greater than $\alpha$. Since $\alpha$ can be taken arbitrarily close to
$d(\nu,\mu,\theta)$, we obtain the lower bound $\dim_H B(\nu) \ge d(\nu,\mu,\theta)$. The inverse
inequality, as was mentioned above, follows from Theorem \ref{7..2}. Thus in the case of inner
point $\nu\in M_1(\Omega)$ the theorem is proved.

If the measure $\nu$ belongs to the boundary of the simplex $M_1(\Omega)$, then one should exclude
from $\Omega$ all elements $i$ with $\nu(i) =0$, and consider the set
\begin{equation*}
 \Omega' =\{\pin i\in \Omega \mid \nu(i)>0\pin\}.
\end{equation*}
Exclude from the brunching process $X_1$, $X_2$, \dots\ all genetic lines containing elements of
colors not in $\Omega'$ and denote as $X'_1$, $X'_2$, \dots\ the resulting brunching process (with
the set of colors $\Omega'$). The corresponding set of infinite genetic lines $X'_\infty$ is
contained in $X_\infty$. It follows from the definition of Billingsley--Kullback entropy
$d(\nu,\mu,\theta)$ that it is the same for the sets of colors $\Omega$ and $\Omega'$. Besides, the
measure $\nu$ lies in the interior of the simplex $M_1(\Omega')$. Therefore, $\dim_H B(\nu)\cap
X'_\infty =d(\nu,\mu,\theta)$ with the same probability as $X'_\infty\ne \emptyset$.

The theorem would be completely proved if the probability of the event $X'_\infty\ne\emptyset$ was
equal to $1-q^*$. But it may be less than $1-q^*$. This obstacle can be overcome as follows.

Let $m' =\sum_{i\in \Omega'} \mu(i)$. This is nothing more than the expectation of each
individual's number of children in the brunching process $X'_1$, $X'_2$, \dots\ If $m'\le 1$ then
\eqref{2,,4} implies the inequality $\rho(\nu,\mu) \ge 0$, which contradicts the condition
$d(\nu,\mu,\theta) >0$ of our theorem. Therefore $m'>1$ and, respectively, the probability of the
event $X'_\infty = \emptyset$ is strictly less than $1$. Let us denote it $q'$.

If the brunching process $X'_1$, $X'_2$, \dots\ was generated not by a unique initial element, but
$k$ initial elements, then the probability of $X'_\infty =\emptyset$ would be equal to~$(q')^k$.
Recall that the cardinality of $X_n$ grows exponentially with probability $1-q^*$. If this is the
case, one can first wait for the event $|X_n|\ge k$, and then consider separately $|X_n|$
independent counterparts of the brunching process $X'_1$, $X'_2$, \dots\ generated by different
elements of $X_n$. This trick allows to obtain the bound $\dim_H B(\nu) \ge d(\nu,\mu,\theta)$ with
conditional probability at least $1-(q')^{k}$ under the condition $|X_n|\to \infty$. Since $k$ is
arbitrary, the above mentioned conditional probability is in fact one, and the complete probability
cannot be less than $1-q^*$. \qed


\end{document}